\input amstex
\input xy
\xyoption{all}
\documentstyle{amsppt}
\document
\magnification=1200
\NoBlackBoxes
\nologo
\hoffset=0.7in
\voffset=0.7in
\def\r{\roman}
\def\p{\partial}
\def\b{\bullet}
\def\D{\Delta}

\def\A{\Cal{A}}
\def\Q{\bold{Q}}
\def\Z{\bold{Z}^+}
\def\R{\bold{R}}
\def\C{\bold{C}}
\def\H{\Cal{H}}

\vsize16cm

%With corrections by Marcolli and Yanofsky

%\hfill{\it Amywork/renorm.tex, version August 23, 2009}

\bigskip

 \centerline{\bf  RENORMALIZATION AND COMPUTATION I:}
 
 \medskip
 
 \centerline{\bf MOTIVATION AND BACKGROUND}

\medskip

\centerline{\bf Yuri I. Manin}

\medskip

\centerline{\it Max--Planck--Institut f\"ur Mathematik, Bonn, Germany,}

\centerline{\it and Northwestern University, Evanston, USA}

\medskip

\hfill{\it To Grigori Mints for his 70th anniversary}

\medskip

\bigskip

{\bf Abstract.} The main observable quantities in Quantum Field Theory,
{\it correlation functions}, are expressed by the celebrated Feynman path integrals
which are not well defined mathematical objects. 
\smallskip
 
{\it Perturbation formalism }
interprets such an integral as a formal series of finite--dimensional but
{\it divergent integrals,} indexed by  Feynman graphs,
the list of which is determined by the Lagrangian of the theory. 

{\it Renormalization}
is a prescription that allows one to systematically  ``subtract infinities''
from these divergent terms producing an asymptotic
series for quantum correlation functions.

On the other hand, graphs treated as  {\it ``flowcharts''},  also form a combinatorial skeleton 
of the abstract computation theory and various operadic formalisms in abstract algebra.

In this role of {\it descriptions} 
of various (classes of) computable functions, such as recursive functions,
functions computable by a Turing machine with oracles etc., graphs can be used to replace
standard formalisms having linguistic flavor, such as  Church's
$\lambda$--calculus and various programming languages. 

The  functions in question are generally not everywhere defined due to
potentially infinite loops and/or necessity to search  in an infinite haystack
for a needle which is not there. 

In this paper I argue that such infinities in {\it classical computation theory}
can be addressed in the same way as 
Feynman divergences, and that  meaningful versions of renormalization in this
context can be devised.  Connections with quantum computation are also
touched upon.

\newpage

\centerline{\bf  Contents}

\medskip

0. Introduction

1. Feynman graphs and perturbation series: a toy model

2. From Feynman graphs to flowcharts and programming methods

3. Bialgebras and Hopf algebras related to computation

4. Regularization and renormalization

 Appendix: Renormalization at large

\bigskip

\centerline{\bf 0. Introduction}

\medskip

{\bf 0.1. Feynman integrals.} The main observable quantities in Quantum Field Theory,
{\it partition and correlation functions}, are expressed by the celebrated Feynman path integrals which are not well defined mathematical objects.  

\smallskip
Perturbation formalism in Quantum Field Theory (QFT)
interprets a  Feynman path integral as a formal series $\sum_{\tau} I_{\tau}$ of
multidimensional integrals, which are generally divergent.
Formal expressions for these divergent terms are 
labeled by (decorated) Feynman graphs $\tau$ (whose exact structure and totality depend
on the Lagrangian of the Quantum Field Theory in question): see  
[Po], [Fra], [Cos], and sec. 2 below based upon   [Ma3], Ch. IV.3.

\medskip

{\bf 0.2. Renormalization.} {\it A renormalization scheme} 
is a prescription (depending on QFT) that allows one to systematically  ``subtract infinities''
from these divergent terms and produce an asymptotic
series for quantum correlation functions. After the initial breakthroughs
made by physicists in the 40s-- 70s,  
 several mathematical versions of renormalization gradually crystallized.
 
 \smallskip

We will use the following  version which
consists of two parts that have essentially different natures.

\smallskip

{\it Part 1: A regularization scheme.}  The divergent integrals corresponding to all individual
Feynman graphs are first ``deformed'', by systematical introduction
of a parameter, say $z$,  such that for (sufficiently small)  $z\ne 0$ 
the integrals converge, whereas for $z=0$ we get an initial
divergent expression. As the simplest example, imagine that
an individual  deformed Feynman integral $I_{\tau}(z)$
defines a germ of a meromorphic function
of $z$ with the only singularity at zero.  Then one can
define  {\it the regularized value of $I_{\tau}$}  as the difference
$$
I_{\tau ,reg}:=(I_{\tau}(z) -\ the\ polar\ part\ of\ I_{\tau}(z))\,|_{z=0}
\eqno(0.1)
$$
(``minimal subtraction'' of counterterms).

\smallskip

The choice of a specific $z$ and the respective deformation 
is an indispensable input of physics at this stage.  
It  has not been formalized in a single concise mathematical scheme:
several known regularizations display remarkable
mathematical and physical variety.

\smallskip

Moreover, even the prescriptions to read off the physically meaningful numbers
from a chosen regularization scheme are wonderfully
ambiguous. 

\smallskip

A number of the standard choices is based on the idea of {\it cut--off}:
e.~g. one introduces a finite scale of integration 
where the infinite one led to divergences (such as ultraviolet
ones), and then calculates, say, an ``effective action'' at this scale.

\smallskip

In this case the infinity of the polar part is not subtracted but rather 
''kept at a safe distance'' 
and  declared to be the bearer of the physical meaning.

\smallskip

A totally non--obvious choice is done in the remarkable schemes of
{\it dimensional regularization}: the dimension of space--time, $4,$
is deformed into a ``complex dimension'' $4+z$, by appropriate
formal changes in integrands. There is a series of other very interesting
examples.

\medskip

{\it Part 2: Renormalization as simultaneous regularization.} 
Feyman integrals corresponding to different
Feynman graphs $\tau$ are mutually interrelated in a way reflecting
the combinatorics of graphs: typically, a divergent integral
contributes to the formal integrand of a larger integral,
that remains divergent even after the (partial) regularization
of the integrand.

\smallskip

This means that the regularization schemes chosen for individual integrals
must be united in a coherent whole.  

\smallskip

After the seminal papers  [Kr1], [ConKr], this can be done in the following way,
in two steps.

\smallskip

First, the linear span of isomorphism classes of Feynman graphs $\tau$
is given the structure of a Hopf algebra $H$. 
Multiplication in this algebra is induced by the disjoint union of graphs,
whereas comultiplication is more sophisticated and roughly speaking,
reflects how bigger graphs are constructed from smaller ones.

\smallskip

Second,
the totality of deformed integrals $I_{\tau}(z)$
induces a character (linear functional)
of the Hopf algebra $H$,   $\varphi :\,H\to A$. In the simplest case
of the regulariztion scheme explained above, this character takes values
in the algebra $A=A_+\oplus A_-$ of germs of analytic functions
meromorphic near $z=0$. Here $A_+$ is the unital
subalgebra of holomorphic germs, and $A_-$ is the subalgebra
of polar parts.

\smallskip

The appropriate characters of $H$ form a group (defined with the help
of the comultiplication of $H$ and the multiplication of $A$). In this group, there is a version of 
Birkhoff decomposition: in particular, $\varphi$ above
can be presented the product of its regular part $\varphi_+$
taking values in $A_+$ and a polar part 
$\varphi_-$ in $1+A_-$. The regular part furnishes all regularized Feynman integrals simultaneously:
$$
I_{\tau, reg}:= \varphi_+(\tau ,z)\,|_{z=0}
\eqno(0.2)
$$

A concise and well--rounded exposition of the relevant mathematics
can be found in [E-FMan].

 \medskip 
 
 {\it  A warning.} The terminology adopted in the sketch above
 is not universally accepted one. The distinction we make between
 regularization and renormalization in physics is generally blurred. 
 In mathematics, the word regularization is used in different senses
 in a variety of contexts having only a vague common 
 intuitive kernel.
 
 \smallskip
 
 Therefore, I summarize once again our current usage.
 
\bigskip

{\it Regularization schemes} in general deal with the problem
of  ``extracting finite information''  from potentially infinite expressions:
summing divergent sums, evaluating divergent integrals, or, as we try to do, taming more general computation
processes. The central, but by no means all--inclusive,  intuitive image is that of
{\it ``subtracting the divergent part''}. The first important examples and insights go back at least to Euler.

\smallskip

{\it Renormalizaton} studies ways of doing it {\it compatibly with composition} of
potentially divergent expressions.

\smallskip

Finally, {\it the Hopf algebra renormalization} adds to this picture the idea that
all relevant compositions can be heaped into one big
Hopf algebra, or group. This
allows one to replace the term--by--term subtraction of divergences
with an overall {\it division} by the divergent part in the respective group
(Birkhoff decomposition).

\medskip

{\bf 0.3. Renormalization and computation.} In this paper 
I argue that certain computational problems could benefit from 
a systematic development of  regularization/renormalization technique.
Moreover, some basic ingredients  of Feynman's renormalization
can be very naturally transposed into  computational contexts.

\smallskip

Here I will start with a brief discussion of such problems, related
Hopf algebras, and regularization schemes.

\medskip

{\bf 0.4. Computation and infinity.} Adopting the classical Church thesis, I
will identify the universe of classical computation theory with that of
{\it partial recursive functions}. Classically, such a function
is a partial map $f:\,\Z\to \Z$ (or, more generally, $f:\,(\Z )^a\to (\Z )^c$).
Values $f(x)$ can be computed in finite running time, if $x$ belongs to the
definition domain $D(f)$ of $f$. However, outside of $D(f)$ the algorithm
computing $f$ (Turing machine, or any other programming method)
will generally work indefinitely long without producing a definite answer.
This is called the undecidability of the Halting Problem.

\smallskip

In order to avoid this kind of  ``computing block'', it is often useful to recede and
replace the initial problem by that of computation of the graph 
$\Gamma_f\subset \Z\times\Z$. The point is
that $\Gamma_f$ is the total image of a {\it general, or total, recursive function},
say $F:\, \Z\to \Z\times\Z$, which is everywhere defined.
Then all points $(x,f(x)), x\in D(f),$ will be printed out in some order.

\smallskip

Both the domain and the target of a partial recursive map $f$ need not be
$\Z$: any enumerable set $X$ having a natural (``computable'') numbering will do.
It is convenient to consider such sets (``constructive worlds'') as objects
of a category (``constructive universe''), and  partial recursive maps 
between them as morphisms.

\smallskip

For example, computability studies in classical analysis usually start
with a working definition of a computable real number.
One can declare, say, that such a number $\alpha$ is represented
by a computable Cauchy sequence computably approximating it, that is, by
a total recursive function $\Z\to \bold{Q},$ $n\mapsto r_n$
such that $|\alpha -r_n| < 2^{-n}.$  A recent paper by M.~Yoshinaga [Yo]
shows that all {\it periods} in the sense of [KoZa] are computable
(see also [Fri] for background). In fact, (certain) periods directly appear
as values of Feynman integrals: cf. [BlEsKr], [Bl], [AlMar1--3] and references
therein.
\smallskip

More generally, challenging problems arise when the target domain is 
(a constructive approximation to) a world of geometric images, say,
subsets of $\R$ or $\R^d$. The case $d=2$ can be motivated by  imagining
computer screen pictures with arbitrary high resolution.   

\medskip

{\bf 0.5. Example: computability of Julia sets.} In [BrYa], the authors
define the constructive world $B^d$ consisting of finite unions of closed balls with rational
centers and rational radii in $\R^d$ (they in fact restrict themselves by dyadic rational numbers 
of the form $p/2^n$.) Such sets  are dense in the Hausdorff metric
on the set of compact subsets. 

\smallskip

A compact set $K\subset \R^d$ is called computable, if there exists
a computable sequence in $B^d$ computably approximating $K$
wrt the Hausdorff metric.  Braverman and Yampolsky then proceed to show,
that for certain {\it computable  numbers} $c$,
the Julia sets $J(z^2+c)$ are {\it not computable.} (Julia sets belong to the class of
most popular fractal pictures). 

\medskip

{\bf 0.6. Example: Code domains and isolated codes.} 
Another class of two--dimensional pictures, {\it code domains,}
naturally arises in the theory of error--correcting codes. The following results
were proved by the author in 1981; for detailed exposition, see cf. [VlMa], Ch. 1, and [TsVlN], sec. 1.3.

\smallskip

Let $F$ be a finite set, {\it alphabet}, of cardinality $q$.
Consider the constructive world of {\it codes}: pairs $(n,C)$ where $n\in \Z$ and $C\subset F^n$. We write also $n(C)$ in place of $n$.
The Hamming distance $d(a,b)$ between two words $a=(a_i),\,b=(b_i)$ in $F^n$ is defined
as the number of $i$ with $a_i\ne b_i$. Define the following
computable functions of codes with integer (resp. rational) values:
$$
d(C)= \roman{min}\,\{d(a,b)\,|\,(a,b)\in C\times C, a\ne b\},\
k(C):= [\r{log}_q\,\r{card}\, C]
$$
$$
R(C):=\frac{k(C)}{n(C)},\quad \delta (C)=\frac{d(C)}{n(C)}.
$$
Denote by $V_q$ the set of all points $(\delta(C), R(C))$ ($q$ being fixed, $n$ variable).
Let $\overline{V}_q$ be the closure of $V_q$ in the square
$[0,1]^2$ of the $(\delta ,R)$--plane. 

\smallskip

It is proved in [VlMa] that   $\overline{V}_q$ consists of all points,
lying below or on some continuous curve $R=\alpha_q(\delta )$,
plus a subset of isolated points lying above this curve.

\smallskip

Similar results were proved for subclasses of codes: {\it linear codes}, which are
linear subspaces of $\bold{F}_q^n$, and polynomial time computable linear codes.

\smallskip

It is not known whether $\alpha_q(\delta )$ is a differentiable function.
Arguably, a high resolution picture could help make an educated guess.
It is not known whether the set of isolated codes or its complement
in ${V}_q$ are computable. 

\medskip

{\bf 0.7. Kolmogorov complexity as ultimate  ``computational infinity''.}
I consider both  problems discussed in 0.5 and 0.6 as  classes of
computations on which the problem of regularization might be tested.
As Braverman and Yampolski ask in [BrYa], ``What would a computer really draw,
when $J(z^2+c)$ is uncomputable?''  Regularization is expected to be
a (theoretical) modification of  computation that produces  better  pictures
at least in most cases. 

\smallskip
I expect that mathematical content of regularization will be 
a  reduction of  other ``computationally challenged'' problems to
the universal one: computation of Kolmogorov complexity. 

\smallskip

The simplest model is this (cf. [Ma2], sec. 5). Denote by $\bold{K}:\,\Z\to\Z$
a (non--computable) permutation rearranging positive integers in the order
of increasing Kolmogorov complexity (with respect to some fixed
optimal numeration). Then any recursive function $f$ becomes bounded by a linear function
in the following sense. There exists a constant $c=c_f$ such that
$$
\bold{K}\circ f \circ \bold{K}^{-1}(n)\le cn\quad \roman{for\ all\ } n\in \bold{K} (D(f)).
$$
Recall that $\bold{K}$ is bounded by a linear function, but $\bold{K}^{-1}$
is not bounded by any recursive function.

\smallskip

Now, imagining the conjugation $\bold{K}\circ f \circ \bold{K}^{-1}$ 
as a regularized version of $f$,
we see that it become ``computable'' in a Pickwickean sense:
regularized values of $f$,
regularized running time and memory of an algorithm, computing  $\varphi$,
all become functions of no more than linear growth.

\smallskip

The following verbal reformulation might serve as a justification of this scheme.
$\bold{K}(x)$ is {\it a short(est) description} of the number $x$, and
$\bold{K}(f(*))$ is a short(est) description of  the number $f(*)$. Thus,
$\bold{K}\circ f\circ \bold{K}^{-1}$ is an avatar of $f$ which replaces 
operations on numbers by operations on their short descriptions. 
An oracle computing $\bold{K}$ and $\bold{K}^{-1}$ would be helpful.

\smallskip

However, regularization in computation remains a challenging problem, and
it is not satisfactorily solved in this paper. (Cf. also  remarks on tropical geometry
in the subsection 4.6).

\smallskip

Much better is the situation with the Hopf algebra analogs.

\medskip

{\bf 0.8. Graphs as programming methods.} Analogs of Feynman 
graphs appear  as  ``flowcharts''  used   for visualization
of information flows in various processes of computation:
cf.  for example [Sc], [Zo], [Mar]. As such, they  often serve
only  illustrative and mnemonic purposes, comparable with the role
of pictures 
in  schoolbook versions of Euclid's geometry.

\smallskip
One notable exception is [Sc] where a class of decorated graphs is structured to form
a Boolean category, and infinite graphs are invoked to encode loops.

\smallskip

Graphs as a category of structured sets in Bourbaki's sense arise
very naturally in the theory of operads of various kinds.
In [BoMa],
we  developed a formalism allowing one to  consider
(decorated) graphs as objects of a monoidal category and to interpret
uniformly many versions of algebraic operads, PROPs etc. as
functors on various categories of graphs.

\smallskip
 As a necessary  step, I show that
the combinatorics involved in the definition of the Hopf
comultiplication for Feynman diagrams can be transported
to flowcharts and in fact to more general categories of {\it programming methods.}

\medskip

{\bf 0.9. Is there an analog of ``action'' in the theory of computation?}
Although the combinatorial skeleton of perturbative series
naturally emerges   in theoretical computing, the latter notoriously lacks
an analog of the basic physical quantity, which is called {\it action.}

\smallskip

Recall that  describing (an isolated piece of) physical reality,
in the classical or quantum mode, we must start with the following two
steps:

\smallskip

a) Define a (mathematical image of) the space  of virtual paths/histories $\Cal{P}$.

\smallskip

b) Introduce a (real--valued) functional $S:\,\Cal{P}\to \bold{R}$ on this space,
satisfying  additivity properties wrt space--time domains (``conservation laws'', ``locality'').

\smallskip

Then in the classical mode of description, physical histories correspond to the
stationary points of $S$ (``principle of the least action''). 

\smallskip
In the quantum mode, calculation of correlators
and transition amplitudes requires understanding of Feynman integrals
$\int_{\Cal{P}}e^{iS(\varphi )} D\varphi$.

\smallskip

The theory of computation badly needs a quantity that
would play the  role of $S$. (Half a century ago, clearly seeing this need, 
I.~M.~Gelfand coined
an expression ``principle of the least interaction'' in the theory
of finite automata).

\smallskip

A version of Kolmogorov complexity again seems to me a promising candidate.

\medskip

{\it  Acknowledgements.}  During the preparation of this paper,
I  profited from correspondence with  Kevin Costello, 
Leonid A. Levin, Matilde Marcolli, Gregori Mints, Alexander Shen, Noson Yanofsky. 
I am very grateful
to them all.

\bigskip

\centerline{\bf 1. Feynman graphs and perturbation series: a toy model}

\medskip

This section is a very brief and elementary introduction to the perturbative
formalism for Feynman path integrals. On a toy finite--dimensional
model we demonstrate, first, how a series over Feynman
diagrams arises, and second, how the structure of its
terms depends on the action functional.

\medskip

{\bf 1.1. Notation.} {\it Feynman path integral} is an heuristic expression of the form
$$
\frac{\int_{\Cal{P}}e^{S(\varphi )}D(\varphi )}{\int_{\Cal{P}}e^{S_0(\varphi )}D(\varphi )}
\eqno(1.1)
$$
or, more generally,  a similar heuristic expression for {\it correlation functions}.

\smallskip

In the expression (1.1), $\Cal{P}$ is imagined as a functional space
of {\it classical fields $\varphi$} on  a {\it space--time manifold} $M$. Space--time may be
endowed with Minkovski or Euclidean metric, but not necessarily:
in models of quantum gravity metric is one of the fields. 
Fields in general may include scalar functions, tensors of various ranks,
sections of vector bundles, connections etc.

\smallskip

$S:\,\Cal{P}\to \bold{C}$ is a functional of {\it classical action}: usually $S(\varphi )$
is expressed as an integral over $M$ of a local density on $M$ which is called {\it Lagrangian.}
In our notation (1.1) (differing by the sign from the standard one),  
$S(\varphi )=-\int_M L(\varphi (x)) dx.$ Lagrangian density may depend on derivatives, include
distributions etc.  

\smallskip

Usually $S(\varphi )$ is represented as the sum of {\it a quadratic part}
$S_0(\varphi )$ (Lagrangian of free fields) and remaining terms which are 
interpreted as interaction and treated perturbatively.

\smallskip
Finally, the integration measure $D(\varphi )$ and the integral itself $\int_{\Cal{P}}$
should be considered as simply a  part of the total expression (1.1) expressing
the idea of  ``summing over trajectories''. 
\smallskip

In our toy model, we will replace $\Cal{P}$ by a finite--dimensional
real space. For the time being, we endow it with a basis indexed by a finite set
of ``colors''  $A$, and an Euclidean metric $g$ encoded by the symmetric
tensor $(g^{ab}),\,a,b\in A.$ We put $(g^{ab})=(g_{ab})^{-1}.$

\smallskip

The action functional $S(\varphi )$ will be a formal series in linear coordinates on $\Cal{P}$, $(\varphi ^a)$,
of the form
$$
S(\varphi )=S_0 (\varphi) + S_1(\varphi ),\quad
S_0(\varphi ):=-\frac{1}{2} \sum_{a,b} g_{ab}\varphi^a\varphi^b,
$$
$$
S_1(\varphi ):=\sum_{k=1}^{\infty}\frac{1}{k!}\sum_{a_1,\dots ,a_k\in A}
C_{a_1,\dots ,a_k}\varphi^{a_1}\dots \varphi^{a_k}
\eqno(1.2)
$$
where $(C_{a_1,\dots ,a_n})$ are certain  symmetric tensors.
If these tensors vanish for  all sufficiently large ranks $n$,
$S(\varphi )$ becomes a polynomial and can be considered as a genuine function on 
$\Cal{P}$.

\smallskip

On the other hand, below we will mostly consider $(g_{ab})$ and  
$(C_{a_1,\dots ,a_n})$ as independent formal variables as well, 
``formal coordinates on the space of theories''.

\smallskip

We can now state our first theorem, expressing the toy version of (1.1)
as a series over (isomorphism classes of)  graphs.

\smallskip
For the time being, a graph $\tau$ for us consists of two finite sets, edges $E_{\tau}$
and vertices $V_{\tau}$, and the incidence map   sending $E_{\tau}$
to the set of unordered pairs of vertices. Each vertex is supposed to be incident to at least one edge. There is one {\it empty graph.}

\medskip

{\bf 1.2. Theorem.}  {\it We have, for a formal parameter $\lambda$
$$
\frac{\int_{\Cal{P}}e^{\lambda^{-1}S(\varphi )}D(\varphi )}{\int_{\Cal{P}}
e^{\lambda^{-1}S_0(\varphi )}D(\varphi )} =
  \sum_{\tau\in\Gamma}\frac{\lambda^{-\chi (\tau )}}{|\roman{Aut}\,\tau |}\,
w(\tau )
\eqno(1.3)
$$
In the r.h.s. of (1.3) the summation is taken over (representatives of) all 
isomorphism classes
of all finite graphs $\tau$. The weight $w(\tau )$ of such a graph
is determined by the action functional (1.2) as follows:
$$
w(\tau ):=\sum_{u:\,F_{\tau}\to A}\ \prod_{e\in E_{\tau}}
g^{u(\partial e)}\prod_{v\in V_{\tau}} C_{u(F_{\tau}(v))}\,.
\eqno(1.4)
$$
Here  $F_{\tau}$ is the set of  flags, or ``half--edges'' of $\tau$.
Each edge $e$ consists of a pair of flags denoted
$\partial{e}$, and each vertex $v$ determines the set  of flags
incident to it denoted $F_{\tau} (v)$.  Finally, $\chi (\tau )$ is the 
Euler characteristic of $\tau$.}

\medskip

{\bf 1.3. Comments.} Even in this toy version, the
meaning of the  ``perturbation series expansion''
(1.3) requires some explanations.

\smallskip

First of all, the numerator of the left hand side of (1.3) is {\it by definition} the result
of {\it term--wise integration} of the formal series which can be obtained
from the Taylor series of the exponent in the integrand. Concretely
$$
\int_{\Cal{P}}e^{\lambda^{-1}S(\varphi )}D(\varphi ) =
\int_{\Cal{P}}e^{\lambda^{-1}S_0(\varphi )} \left(1+\sum_{N=1}^{\infty}
\frac{\lambda^{-N}S_1(\varphi )^N}{N!}\right)\, \prod_a d\varphi^a \  :=
$$
$$
\int_{\Cal{P}}e^{\lambda^{-1}S_0(\varphi )} \prod_a d\varphi^a \ +
$$
$$
\sum_{N=1}^{\infty}\frac{\lambda^{-N}}{N!}\sum_{k_1,\dots ,k_N=1}^{\infty}\frac{1}{k_1!\dots k_N!}
\sum_{a^{(i)}_j\in A, 1\le j\le k_i}\prod_{i=1}^N
C_{a_1^{(i)},\dots ,a_{k_i}^{(i)}}
\int_{\Cal{P}} e^{\lambda^{-1}S_0(\varphi )}\prod_{i,j}^N \varphi^{a_j^{(i)}}
 \prod_a d\varphi^a  \,.
  \eqno(1.5)
$$
\medskip

This definition makes sense if the right hand side of (1.5) is understood
as a formal series of infinitely many independent weighted variables $C_{a_1,...,a_k}$,
weight of  $C_{a_1,...,a_k}$ being $k$. In fact, the Gaussian
integrals in the coefficients uniformly converge.

\smallskip

More precisely, putting $d:=\roman{dim}\,\Cal{P}=\roman{card}\,A$ and 
$D:= \roman{det}\,(g_{ab})$, we have
$$
\int_{\Cal{P}}e^{\lambda^{-1}S_0(\varphi )} \prod_a d\varphi^a \ =
\frac{(2\pi\lambda )^{d/2}}{D^{1/2}}.
$$
Furthermore, put for a polynomial $F(\varphi )\in \bold{C}[\varphi ]$
$$
\langle F(\varphi )\rangle :=\frac{\int_{\Cal{P}}e^{\lambda^{-1}S_0(\varphi )} 
F(\varphi ) \prod_a d\varphi^a}{\int_{\Cal{P}}e^{\lambda^{-1}S_0(\varphi )} \prod_a d\varphi^a}\ .
$$
Then we have the following

\medskip

{\bf 1.3.1. Wick's Lemma.} {\it In the notations above, we have
\smallskip

a) $\langle \varphi^{a_1}\dots \varphi^{a_n}\rangle =0$ for $n \equiv 1\,
\roman{mod}\,2\,.$
\smallskip
b) $\langle \varphi^{a}\varphi^b\rangle = \lambda g^{ab}.$
\smallskip
c) $\langle \varphi^{a_1}\dots \varphi^{a_{2m}}\rangle =
\lambda^m\sum g^{a_{i_1}a_{j_1}}\dots g^{a_{i_m}a_{j_m}}$
where the summation is taken over all unordered partitions of
$\{1,\dots ,2m\}$ into $m$ unordered pairs $\{i_1,j_1\},\dots ,
\{i_m,j_m\}$ (pairings.)}

\medskip

{\bf 1.4. Proof of Theorem 1.2.} Let us now calculate the l.h.s. of (1.3). From the definition and (1.5) we get 
$$
\sum_{N=1}^{\infty}\frac{\lambda^{-N}}{N!}\sum_{k_1,\dots ,k_N=1}^{\infty}\frac{1}{k_1!\dots k_N!}
\sum_{a^{(i)}_j\in A, 1\le j\le k_i}\prod_{i=1}^N
C_{a_1^{(i)},\dots ,a_{k_i}^{(i)}}\,\,
\langle \prod_{i,j}^N \varphi^{a_j^{(i)}}\rangle \,.
  \eqno(1.6)
$$

Choose some $(N; k_1,\dots ,k_N)$. A typical monomial of degree $N$ in
$C_{\bullet}$ in the decomposition
of (1.6) will be
$$
\lambda^{-N}\,\frac{1}{N!}
\prod_{i=1}^N\frac{1}{k_i!}\,
 C_{a_1^{(i)},\dots ,a_{k_i}^{(i)}}\,\,
\langle \prod_{i=1}^N\varphi^{a_1^{(i)}}\dots\varphi^{a_{k_i}^{(i)}}
\rangle \,.
\eqno(1.7)
$$
It vanishes if $k_1+\dots +k_N$ is odd. Otherwise, in view of
Wick's Lemma (1.7) can be rewritten as
$$
\lambda^{-N+\frac{1}{2}\sum k_i}\,\frac{1}{N!}
\prod_{i=1}^N\frac{1}{k_i!}\,
 C_{a_1^{(i)},\dots ,a_{k_i}^{(i)}}
\left( \sum g^{a_{l_1}^{(i_1)}a_{m_1}^{(j_1)}}\dots
g^{a_{l_r}^{(i_r)}a_{m_r}^{(j_r)}}
\right) 
\eqno(1.8)
$$
where $r=\frac{1}{2} \sum k_i$ and the inner sum is taken over all pairings
of the set of ordered pairs $F=F(N;\,k_1,\dots ,k_N)=\cup_{i=1}^N\{ (i,1),\dots ,(i,k_i)\}.$

\smallskip

Construct now a family of graphs $\tau$, corresponding to  the
monomials in $ g^{**}$ appearing  in (1.8). They will have
a common set of flags $F_{\tau}:=F,$ and a common set of  vertices 
$V_{\tau}=\{1,\dots ,N\}$, $\partial_{\tau}(i,l)=i$.  Declare that
two flags constitute halves of an edge, if these flags are paired
as in the respective monomial  in (1.8).
If we color the flags of one such graph by the map $F_{\tau}\to A:
(i,l)\mapsto a^{(i)}_l,$ then we will get a monomial
in  $(C_{\bullet}, g^{**})$ appearing in the weight function (1.4).

\smallskip
It remains
to perform some bookkeeping in order to identify
the coefficients at this monomials appearing respectively
in (1.6) and the r.h.s. of (1.3).  Here is a sketch.

\smallskip

The graphs constructed above  represent all isomorphism classes
of graphs in our sense.  In fact, a choice of $(N;\,k_1,\dots ,k_N)$
determines the number of vertices of any valence,
and the choice of a pairing determines which pairs of
flags become edges ($N=0$ produces the empty graph.)
Moreover, a non--empty graph comes thus equipped with a total ordering
 of its vertices
and all sets of flags belonging to one vertex.
The sum over graphs does not take care of these orderings.
The group $\roman{Aut}\,\tau$ effectively acts on the whole set
of them consisting of $N!\,\prod_{i=1}^Nk_i!$ elements.
Summing over isomorphism classes, we may replace the numerical
coefficient in (1.8) by $|\roman{Aut}\,\tau |^{-1}.$

\smallskip

Finally, 
$$
-N+\frac{1}{2}\sum_{i=1}^Nk_i=-|V_{\tau}|+|E_{\tau}|=\chi (\tau ).
$$
\smallskip

This completes the proof.

\medskip

  There are several more useful identities
in the framework of our toy model. We mention two of them;
proofs can be found in [Ma3], IV.3.

\medskip

{\bf 1.5.1. Theorem.} {\it We have
$$
\roman{log}\,\, \sum_{\tau\in\Gamma}\frac{\lambda^{-\chi (\tau )}}{|\roman{Aut}\,\tau |}\,
w(\tau )=
 \sum_{\tau\in\Gamma_0}\frac{\lambda^{-\chi (\tau )}}{|\roman{Aut}\,\tau |}\,
w(\tau )
\eqno(1.9)
$$
where $\Gamma_0$ is a set of representatives of isomorphism classes of connected graphs.}
\medskip

\medskip

{\bf 1.6. Summation over trees and the stationary value of the classical action.}
It is well known that quantum path integral is expected to be dominated
by small fluctuations around stationary points of the classical
action functional $S(\varphi )$. In our formal context, one can define  a natural
stationary point $\varphi_0$, and it turns out that $S(\varphi_0)$ is again
a sum over graphs with non--empty set of edges, this time simply connected and connected
that is, trees $T$.

\smallskip

We will again treat $C_{a_1,\dots ,a_k}$ as independent graded formal variables
over a ring containing $g_{ab},g^{ab}$ and $\bold{Q}.$
Then all our sums make sense as formal series.

\smallskip

Put $C^a=\sum_{b\in A}g^{ab}C_b$ and denote by $N\subset R$
the ideal generated by $C_{a_1,\dots ,a_k}$ for all $k\ge 2.$

\medskip

{\bf 1.7. Theorem.} {\it a) The equations
$$
\forall a\in A,\ \frac{\partial S (\varphi )}{\partial \varphi^a}=0
\eqno(1.10)
$$
admit the unique solution $\varphi_0=\{\varphi^a_0\}$, $a\in A$,
satisfying the condition
$$
\varphi_0^a \equiv C^a\,\roman{mod}\, N\,.
\eqno(1.11)
$$
\smallskip

b) The series over trees 
$$
Z:= \sum_{\tau\in T}\frac{\lambda^{-\chi (\tau )}}{|\roman{Aut}\,\tau |}\,
w(\tau )
\eqno(1.12)
$$
satisfies the differential equations
$$
\frac{\partial Z}{\partial C_a}=\varphi_0^a,\ a\in A,
\eqno(1.13)
$$
and is the respective critical value of $S (\varphi ):$}
$$
Z=S (\varphi_0)\,.
\eqno(1.14)
$$
\medskip

{\bf 1.8. Quantum fields: general indications.}  Using the toy model as an inspiration for developing
the respective perturbative formalism for more realistic
quantum field theories, one meets many new problems. The simplest of them
can be illustrated already in the case of scalar field theories.

\smallskip
 
Basically, in this case our combinatorial finite set of indices $A$ is replaced
by the set of points of space--time $M$, and a vector
$\varphi =\{\varphi^a\}$ by a function  $\varphi =\varphi (x)$, $x\in M$.
The right hand side of (1.3), suitably interpreted, becomes 
{\it a definition} of the path integral.
When one devises this interpretation, summations over $A$ 
are replaced by integrations over $M$.
In particular, weights of Feyman's graphs (1.4) turn into multidimensional
integrals, and these integrals generally turn out to be  {\it divergent.}

\smallskip

The simplest divergences occur already in the interpretation  of (the products of) the terms
$g^{ab}$. When $S_0(\varphi )$ is given by the Lagrangian
density of the free field  on $\bold{R}^D$
$$
L(\varphi ):=\frac{1}{2}(\sum_{k=1}^d|\partial_k\varphi (x)|^2 +m^2\varphi (x) ^2),
$$
the classical ``stationary phase'' equation is Klein--Gordon
$(-\Delta+m^2)\varphi (x)=0$, and passing to the Fourier transform in the
momentum space, we easily get an integral representation for the relevant continuous
analog of $(g^{ab})$, ``free propagator'':
$$
G_0(x-y)=\int_{\bold{R}^D}\frac{e^{-ip\cdot (x-y)}}{p^2+m^2} \frac{d^Dp}{(2\pi )^D}
$$
For $D=1$, it is a continuous function, but for $D>1$, it is a distribution singular on the diagonal
$x=y$. Trying to make sense of the product of terms $g^{ab}$ corresponding, for example,
to the graph with two vertices, connected by two edges, and supplied with one aditional flag 
at each vertex, we get a formal expression 
$$
\int d^Du\, d^Dv\, G_0(x-y)\, G_0(u-v)^2\, G_0(v-z)
$$
which after  Fourier  transform in $q$, with fixed $p$, becomes 
$$
\int \frac{1}{q^2+m^2} \frac{1}{(p-q)^2+m^2}  \frac{d^Dq}{(2\pi )^D}.
\eqno(1.15)
$$
This integral diverges for $D\ge 4.$

\smallskip

Therefore, such integrals must be suitably regularized. 

\smallskip

The choice of a regularization method for individual
integrals (weights $w(\tau )$), such as introducing an auxiliary parameter and then
defining a minimal subtraction of divergent terms, must take into account
the basic fact that the perturbation series contains
{\it all graphs}, and if one say, subtracts the divergence
of weight of some subgraph, this changes the weight
of the total graph and its initial divergence.

\smallskip

Beautiful inductive combinatorial prescriptions for doing subtractions
coherently  were devised by physicists
in the early years of path integration.  Nowadays they are organized into a
well--structured algebraic theory involving Hopf algebras.

\smallskip

Basically, comultiplication in this Hopf algebra reflects
the ways a graph can be composed from its subgraphs, and weight
function is translated into a character of this Hopf algebra.
Subtraction of divergences is replaced by dividing out the polar part
in the group of characters. The following sections of the article
are dedicated to various implementations of this general scheme.

 \smallskip
 
 In conclusion, we remark that integrals of the type (1.15), when they absolutely converge,
 as well as the regularized values of more general Feynman integrals,
 are {\it periods} in the sense of Kontsevich--Zagier ([KoZa]), and therefore
 are closely related to (mixed) motives. Study of various classes of such integrals
 from this perspective is very active: see, in particular, [AlMar1--3], [Bl],
 [BlEsKr], and references therein.

\bigskip

\centerline{\bf 2. From Feynman graphs to flowcharts and programming methods}

\medskip

 In this section, we review several contexts in which  graphs
 and their weights such as (1.4) appear without reference to path integrals.
 We stress their interpretation as flowcharts processing
 input data into output data, and briefly discuss two generalizations:
{\it programming methods}  and (generalized) {\it operads.} 
 This requires introduction of more general
 notion of decorated graphs. We start with some preparatory remarks.

 \medskip
 
 {\bf 2.1. Weights of connected graphs and tensor networks.}
 The expression (1.4) for a toy Feynman weight
$$
w(\tau ):=\sum_{u:\,F_{\tau}\to A}\ \prod_{e\in E_{\tau}}
g^{u(\partial e)}\prod_{v\in V_{\tau}} C_{u(F_{\tau}(v))}\,.
\eqno(2.1)
$$
can be slightly generalized and reinterpreted as a polylinear map
$$
\otimes_{v\in V_{\tau}} S^{|F_{\tau}(v)|} (\Cal{P}) \to k.
\eqno(2.2)
$$
Namely, let the former finite--dimensional ``path'' vector space
$\Cal{P}$ be now defined over a field $k$ and endowed with symmetric
non--degenerate metric $g:\,S^2(\Cal{P})\to k.$ Let
$C_v\in S^{|F_{\tau}(v)|} (\Cal{P})$ be a family of symmetric tensors
indexed by vertices of $\tau$. 

\smallskip
Such a family can be interpreted
as a {\it tensor network}: cf. [Zo]. The weight (2.1)
is then the result of contraction of all these tensors
along all pairs of indices corresponding to (pairs of flags forming)
edges of $\tau$. More precisely, factors $g^{u(\partial e)}$ raise subscripts, and summation over $u$ comes from writing tensors in a basis. 
 
\medskip

{\bf 2.2. Variants and generalizations.} (i) First,  expressions (2.1)
and (2.2) suggest to extend the notion of graph by allowing
``freely hanging'' flags (they are variously called leaves, tails, or legs).
Then the basic constituents of any graph are ``corollas''.

\smallskip
By definition, a corolla is an one--vertex graph 
 with several flags attached to it  at one end. Supplying
 an additional information specifying which flags are paired
 to form halves of an edge, we can construct any graph
 from the corollas of all its vertices.
 
 \smallskip
 
 (ii) Second, graphs should be considered as objects of a category. 
 We have already used a notion of isomorphism/automorphism
 of graphs which clearly involved only the combinatorial part of the structure.
 
 \smallskip
 Developing this idea further, we can introduce other classes of morphisms.
 One may, for example, interpret as the morphism a  list of pairs of tails
 of a given graph (source) that must be glued together in order to produce 
  edges out of a new graph (target). Another class of
 morphisms consists of contractions of a subset of edges. 
 
 \smallskip
 
 Moreover, one is naturally led to the consideration of {\it decorated
 graphs}, for example, oriented ones. This is essential when graphs are
 interpreted as flowcharts. In fact, in (2.2) we could bypass
 orientation only because the metric made the space $\Cal{P}$
 self--dual. 
 
 \smallskip
 In more general models, the set of flags is subdivided
 into two subsets: inputs (oriented towards the vertex) and outputs (oriented outwards).
 If an input $f$ is decorated by a linear space, the output $f^{\prime}$
 forming the other half of the edge must be decorated by the dual space,
 so that ``contraction of indices'' will still be possible. 
 
 \smallskip
 
 Then morphisms must be compatible with decorations.
 The compatibility conditions are usually motivated by envisioned 
 applications. 
 
 \smallskip
 
 (iii) Third, having set upon a category of decorated graphs, one  naturally 
 wonders what are interesting functors defined on graphs.
 It turns out that all known types of {\it operads}
 can be treated as such functors, and operadic algebras
 become natural transformations of functors. 
 In particular, graphs of the type that appeared in the toy model
 produce {\it cyclic operads}, and tensors $\{g^{ab},C_{a_1,\dots ,a_n}\}$
can be treated as structure constants of respective algebras.

\smallskip

A systematic development of this formalism is given in [BoMa].
See also [vdL], [Va1--2], [I1--2], [Bl] for various other versions.
Pseudo--tensor categories of [BeDr] can be introduced 
in a similar way.

\smallskip

Below, I will sketch only those basic definitions that are relevant in the
computation/renormalization context.
 
\medskip

{\bf 2.3. Combinatorial graphs and geometric graphs.}
We make a distinction between {\it combinatorial graphs}
which are certain Bourbaki structures defined via discrete (in our case even finite)
sets, and their {\it geometric realizations}, which are
topological spaces. Combinatorial graphs form a category $Gr$.
A general description of morphisms of $Gr$ is given in [BoMa].
Some  subclasses of morphisms
will be described where they are needed.

\smallskip

{\it A combinatorial graph $\tau$} consists of two sets
$F_{\tau}$ (flags), $V_{\tau}$ (vertices)  and two incidence relations.
Each flag  $f\in F_{\tau}$ is incident to exactly one vertex $v\in V_{\tau}$,
its {\it boundary} which is $v=\partial_{\tau}(f),$ and 
the map $\partial_{\tau}:\,F_{\tau}\to V_{\tau}$
is a part of the data.  Finally, some pairs of flags form ``halves of an edge'': this incidence relation is represented
by an involution $j_{\tau}:\,F_{\tau}\to F_{\tau}$, $j_{\tau}^2=\roman{id}.$

\smallskip

$$
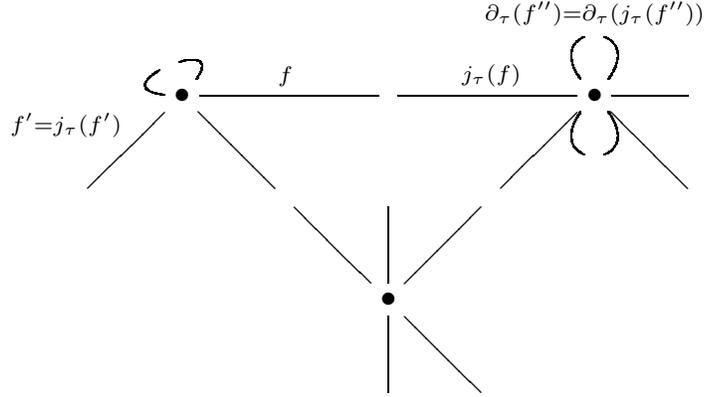
\xymatrix{
& \bullet\ar@{-}[rr]^f
\ar@{-}[dr]\ar@{-}[dl]_{f'=j_{\tau}(f')}\ar@{-}@(l,ur)[]|\hole&& \ar@{-}[rr]{ }\ar@{-} [rr]^{j_{\tau}(f)} && \bullet\ar@{-}[r]\ar@{-}[dr] \ar@{-}[dl]
\ar@{-}@(ur,ul)[]|\hole_{\partial_\tau (f'')=\partial_\tau (j_\tau (f''))} \ar@{-}@(dr,dl)[]|\hole&
{ }
\\
{ }&& { } &{ }& { } && { }
\\
&&&\bullet\ar@{-}[d]\ar@{-}[rd]\ar@{-}[ur]\ar@{-}[ul]\ar@{-}[u]\\
&&& { }&{ }
}
$$

{\it A geometric graph}, geometric realization  of $\tau$, is the topological space
$|\tau |$ which is obtained from the family of  segments $\{[0,1/2]_f\}$ 
indexed by $f\in F_{\tau}$.  We have to identify
some groups of end--points. Namely, we  glue together $0$'s of all segments,
whose respective flags have one and the same boundary vertex $v$. Finally, we
glue the end--points $1/2$ of each
pair of segments indexed by two flags $f\ne f^{\prime}$ such that 
$j_{\tau}(f)=f^{\prime}.$

\smallskip

Pairs of flags  $f \ne f^{\prime}$ with $j_{\tau}(f)=f^{\prime}$ 
are elements of the set of {\it edges} $E_{\tau}$. Fixed points
of $j_{\tau}$ are called {\it tails}, their set is denoted $T_{\tau}.$

\smallskip
A graph $\tau$ is called {\it connected} (resp. {\it simply connected}, etc) iff its geometric realization is such.  A connected and simply connected graph is a {\it tree},
a disjoint union of trees is a {\it forest.} A tree with one vertex is a {\it corolla.}
Each vertex $v\in V_{\tau}$ defines its corolla $F_{\tau}(v)$ consisting
of flags incident to $v$, $v$ itself, induced $\partial_{\tau}$ and trivial $j_{\tau}.$

\smallskip

We will not consider isolated vertices, so that our $\partial_{\tau}$
will be a surjection. However, we do allow  {\it empty graphs.} 

\smallskip

Finally, a few words about {\it morphisms.} Any morphism
of combinatorial graphs $h:\,\tau \to \sigma$ that we will consider, 
will be {\it uniquely defined by two maps}:
$$
h_V:\, V_{\tau} \to  V_{\sigma},\quad  h^F:\, F_{\sigma} \to  F_{\tau}
$$
However, conditions, restricting allowed maps, will depend on the class of morphisms,
and eventually on the decorations (cf. below).
Composition is composition of maps.

\smallskip

In particular, $h$ is {\it an isomorphism}, iff $h_V$, $h^F$ are {\it bijections,
identifying the incidence maps.}

\smallskip

Notice  one peculiarity: $h_V$ is covariant, whereas $h^F$ is {\it contravariant.}
This choice can be explained using the  flowcharts intuition discussed below:
 a change of arguments produces the lift of functions
{\it in reverse direction.}

\medskip

{\bf 2.4. Decorated graphs.}  Let $L=(L_F, L_V)$
be two sets: {\it labels of flags and vertices}, respectfully.

\smallskip

{\it An $L$--decoration} of the combinatorial graph $\tau$ consists of
two maps $F_{\tau}\to L_F$, $V_{\tau} \to L_V$. Usually these maps
are restricted by certain {\it compatibility with incidence relations} 
conditions that we will not try to
axiomatize, and only illustrate on several basic examples.

\smallskip

$L$--decorated graphs for various specific sets of labels 
will also form a category. Any morphism $h$ will be, as above,
determined by $h_V, h^F$, but this time these maps will
be restricted by certain {\it compatibility with decorations.} 

\smallskip

{\it An isomorphism} of decorated graphs is simply an isomophism
of underlying combinatorial graphs,  preserving the decorations.

\medskip

{\bf 2.5. Orientation and flowcharts.} Consider the
 set of labels $L_F:=\{in, out\}$. A decoration
$F_{\tau} \to \{in, out\}$ such that
halves of any edge are oriented by different labels, is called {\it an
orientation of $\tau$.} On the geometric realization,
a flag marked by $in$ (resp. $out$) is oriented towards
(resp. outwards) its vertex.

\smallskip

Tails of $\tau$ oriented $in$ (resp. $out$) are called
{\it (global) inputs $T_{\tau}^{in}$} (resp. {\it (global) outputs} $T_{\tau}^{out}$) of $\tau$.
Similarly, $F_{\tau}(v)$ is partitioned into inputs and outputs
of the vertex $v$.

\smallskip

Consider an orientation of $\tau$. Its edge is
called {\it an oriented loop}, if both its halves belong
to the same vertex. Otherwise an oriented edge
starts at a source vertex and ends at a different
target vertex.

\smallskip

More generally, a sequence of pairwise distinct
edges $e_1, \dots ,e_n$, is called {\it a simple path}
of length $n$, if $e_i$ and $e_{i+1}$ have a common vertex, and
the $n-1$ vertices obtained in this way are
pairwise distinct. If moreover $e_1$ and $e_n$
also have a common vertex distinct from the mentioned ones,
this path is {\it a wheel} of length $n$.
A loop is a wheel of length one.

\smallskip

A wheel in an oriented graph is called {\it oriented wheel},
if  of each two flags of a wheel, sharing a common vertex, one is oriented $in$, and 
and another $out$.
\smallskip

Imagine now that we have a set of operations $Op$ that can be performed on certain
objects from a set $Ob$ and produce other objects of this kind. Take
(names of) this set of operations 
as the set of labels of vertices $L_V$. Then an oriented graph $\tau$
with vertices decorated by  $L_V$ can sometimes be considered
as a flowchart describing a set of computations.

\smallskip

Here is the simplest case.  Let $\tau$ be an oriented corolla,
with vertex labeled by the name of an operation $f$
which describes a map $Ob^a\to Ob^c$, where $a,c\in \bold{N}$.
Call $a $ (resp. $c$) the {\it arity} (resp. {\it co--arity}) of $f$. 
Require the following compatibility: 
$$
\roman{card}\,T_{\tau}^{in}=a,\quad   \roman{card}\,T_{\tau}^{out}=c.
$$
Enrich the decoration of  inputs by a bijection   $T_{\tau}^{in}\to \{1,\dots ,a\}$
and that of outputs by a bijection $T_{\tau}^{in}\to \{1,\dots ,c\}$.
Then such a corolla can be considered as a replacement
of the expression $f(x_1,\dots ,x_a)=(y_1,\dots ,y_c)$ with variable $x_i, y_j\in Ob$.

\smallskip

More general flowcharts $\tau$ are obtained, if we allow outputs of a set of operations
to be taken as inputs of another set of operations. At the end we will  get
global outputs.

\smallskip

Notice that if we choose $Ob \times \{in,out\}$ as the label set for decorating
flags $L_F$,  then certain decorated graphs can be interpreted
as concrete {\it instances} of a computation.

\smallskip

Three remarks are in order here.

\smallskip

(i) An oriented corolla with  {\it empty}  set of inputs, but non--empty set
of outputs, can be used to formalize the notion of {\it oracle} in our context.
In theoretical computer science, one is allowed to imagine an oracle
providing one with a piece of information which did not come
from any sensible computation. 

\smallskip

(ii) One can also imagine a corolla with empty set of global outputs. 
But a more intuitive formalization of such a ``device'' is a corolla with one output
which always produces a specific ``empty''  object. Its vertex may also 
be decorated by a name of ``empty'' operation. 

\smallskip

(iii) If our purported flowchart has oriented loops or wheels, we may be in trouble, 
because intuitively they require explicit inclusion of ``time'' in our calculation:
the vertex of a loop must feed its output into its input -- how many times?
When this circling of information will stop?

\smallskip

Without adopting specific prescriptions, we better avoid oriented wheels
in our flowcharts by using only directed graphs (but cf. the definitions in [Sc]).

\smallskip

An oriented graph $\tau$ is called {\it directed} if it
satisfies the following condition:
\smallskip

($\bullet$) {\it On each connected component of the geometric realization $|\tau |$,
one can define a continuous real valued function (``height'')
in such a way that
moving in the direction of orientation along each flag
decreases the value of this function.}

\smallskip

In particular, oriented trees and forests are always directed.

\smallskip

Of course, a directed graph admits infinitely many compatible height functions,
but only a partial order induced by such a function on vertices,
will be used below.  However, if difference of heights 
models time of computation, then the function itself becomes 
an essential element of the structure.

\medskip

{\bf 2.6. Example: Dana Scott's flow diagrams.} In [Sc], a class of
(eventually infinite) decorated graphs is introduced called {\it flow diagrams.}
It was explicitly designed to model computations with loops,
possibly infinite ones, and it forms a special kind of category which
Dana Scott treats as a lattice.

\smallskip

Here we give only basic definitions.

 \medskip

{\bf 2.6.1. Labels.} The set $L_V$ (labels of vertices) has the following structure:
$$
L_V=\Cal{F}\coprod  \Cal{S}\coprod \{\bot , \top\}.
$$
Here $\Cal{F}$ are functional labels: names of {\it operators} transforming
its input to its output. This set includes the name $I$
of identical operator. Furthermore, $\Cal{S}$ are names
of  {\it switches}: a switch tests its input and redirects it
to {\it one of its outputs,} depending on the results.
Finally, symbols $\bot$ and $\top$ are used to generate some
``improper'' diagrams. In particular, a vertex labeled by $\bot$
describes ``the smallest'' flow diagram, whereas $\top$
corresponds to the ``largest'' one. 
\medskip

{\bf 2.6.2. Graphs.} The first subclass of decorated graphs qualifying as 
Scott's flow diagrams are {\it oriented trees}. The orientation
(describing the direction of the information flow) goes from
one tail (input) to many tails (outputs). Each vertex has one input
and either one output, or two outputs. If it has one output, it must be
labeled by a symbol from $\Cal{F}$ or else $\top ,\bot$.
If it has two outputs, it must be labeled by a switch from $\Cal{S}$.

\smallskip

Clearly, when such a finite tree processes a concrete input, the output will show
{\it only at one of the output tails}, because of the semantics of switching.
Hence we might imagine a special vertex accepting many inputs 
and producing the one which is ``not silent''. A price for it will be that our graphs will 
not be trees anymore. They will still remain {\it directed graphs.}

\medskip

We now pass to another example.

\medskip

{\bf 2.7. Example: Yanofsky's algorithms.} Recursive functions
are certain (partial) maps $f:\,(\Z)^a\to (\Z)^c$. The number of arguments 
$a$ is called  {\it arity} of $f$, $c$ its {\it coarity}. The definition domain
of $f$ is denoted $D(f)$. 

\smallskip

{\it Basic recursive functions} are: successor $x\mapsto x+1$,
projections, and constant functions. The world of recursive functions
is the smallest set of functions containing all basic
functions and closed with respect to four {\it elementary operations}: composition
(or substitution), bracketing, recursion
and the $\mu$--operator. The first three, applied to everywhere defined recursive functions, produce an everywhere defined function as well.
The $\mu$--operator generally creates only a partial function.
For more details, see e.~g.  [Ma1], V.2, and below. 

\smallskip
 Recursive functions can be  introduced by their {\it descriptions}
 in a formal or programming language: essentially, 
 a sequence of functions whose first term is (the name of) a basic
 function, and the name of a new function is supplied by the name
 and arguments of elementary operation, applied to some formerly constructed functions.
  If the operator $\mu$
 is not used in a description, the resulting function is called
 {\it primitive recursive}.

\smallskip

In [Ya], N. Yanofsky suggested  using graphs (modulo
appropriate equivalence relation) as a replacement of the 
descriptions above, at least for primitive recursive functions.
I will present  some of his constructions below to illustrate
our general approach.
 
 \medskip

{\bf 2.8. Decorated graphs $Prim$.}
Elements of $Prim$ 
are disjoint unions of trees $\tau$, in which each vertex is the boundary of at least two flags. Moreover, $\tau$ must be endowed
with an {\it admissible decoration.} The latter consists of the following data. They can be chosen independently  on each connected component so that in the following discussion
we speak about trees, if we did not explicitly mentioned the general case.

\smallskip

(a) {\it A marked tail}, which is called {\it root},
or {\it the (global) output} of $\tau$. Its vertex is called
the {\it root vertex}. The remaining tails are called {\it (global)
inputs} of $\tau$. Global inputs form a set
$F_{\tau}^{in} \subset F_{\tau}$, and we consider the global output
as an one--element subset $F_{\tau}^{out} \subset F_{\tau}$.

\smallskip

A choice of root determines (and is equivalent to) the choice of  a
specific {\it orientation}: a map $F_{\tau}\to \{in,out\}$.
Namely, in each shortest path (sequence of flags) 
from a global input
to the root, assign $out$ to the flag that leaves its vertex,
and $in$ to the flag that enters it. This defines the partition
of all flags into two subsets: (local) inputs and outputs.

\smallskip

We will say that $\tau$ with such a decoration is {\it an oriented tree.} We repeat that by definition, each oriented tree
must have exactly one global output, and at least one global
input. 
\smallskip

(b) All corollas of an oriented tree are also oriented trees.
The next part of decoration is a choice of
{\it total order on the set of inputs of each corolla of $\tau$,}
and, if $\tau$ is not connected, a choice of total order on the set of its connected component.  

\smallskip

(c) A map {\it arity/coarity}: $F_{\tau}\to \bold{N}:\,f\mapsto (a(f),
c(f))$. 
If two flags are halves of an edge, they must be assigned the same arity/coarity.

\smallskip

(d) A map $op:\, V_{\tau}\to \{\bold{c},\bold{b},\bold{r}\}$.
The value $op\,(v)$ assigned to a vertex is called the respective {\it operator}: $\bold{c},\bold{b},\bold{r}$
stand respectively for {\it composition, bracketing, recursion}.

\smallskip

(e) A map $in:\, F_{\tau}^{in}\to\ \{basic\ recursive\ functions\}$
such that for each $i\in F_{\tau}^{in}$, $in(i)$ is a basic function
of arity $a(i)$ and coarity $c(i)$.

\smallskip

All these data must be {\it compatible}. A part of compatibility conditions
was already included in the description.
We will now formally introduce the remaining set,
and simultaneously explain an interpretation of
graphs in $Prim$ (without decoration 2.8 (e)) as operations 
acting on families of input functions.

\medskip

{\bf 2.9. Objects of $Prim$ as flowcharts.} 
Given an oriented tree $\tau$ with a  decoration as above,
we interpret the whole $\tau$ as a symbol of an {\it operation
$Op (\tau )$}
that can be performed over  families
of functions, indexed by
global inputs of $\tau$.     

\smallskip

More precisely, let $f=\{f_i\,|\,i\in F_{\tau}^{in}\}$
be a family of  functions (or even partial functions)
such that $f_i:\, (\bold{Z}^+)^{a(i)}\to (\bold{Z}^+)^{c(i)}.$ 
Then 
$$
Op\, (\tau )(f)= g:\ (\bold{Z}^+)^a\to (\bold{Z}^+)^c
$$
where $(a,c)$ is the arity/coarity of the root.
\smallskip

The prescription for getting $g$, given $f$, runs as follows.

\medskip

{\it One vertex case.} Let $\tau$ be a corolla whose vertex is decorated by 
$\bold{c},\bold{b},$ or $\bold{r}$. Then $g$ is obtained
by applying to the family $\{f_i\},\, i\in F_{\tau}^{in}$, the respective elementary operation:
composition, bracketing or recursion.
This requires the following compatibilities which vary depending
on the label of the vertex.

\medskip

(a) {\it Composition.}  
Let $(a_1,c_1),\dots ,(a_r,c_r)$ be the family of arities/coarities
of inputs ordered as the respective flags.
They must then be constrained by the condition
$c_1=a_2,\dots ,c_{r-1}=a_r$, and the arity/coarity 
of the output must be $(a_1,c_r)$.

\smallskip

For a general $\tau$, these compatibility conditions must be
satisfied for all corollas $\tau_v$ of all vertices decorated by 
$\bold{c}$.

\smallskip

In the flowchart interpretation, such a corolla 
transforms an input family $(f_1,\dots,f_r)$,
$f_i:\,(\bold{Z}^+)^{a_i} \to (\bold{Z}^+)^{c_i}$, into
the composition $f_r\circ f_{r-1}\circ \dots f_1$.

\medskip

(b) {\it Bracket.}  
With the same notation as  in (a),
the compatibility condition reads $a_{\bullet}:=a_1=\dots =a_r$,
and the arity/coarity of the output must be
$(a_{\bullet},c_1+\dots c_r)$.

\smallskip

For a general $\tau$, these compatibility conditions must be
satisfied for all corollas $\tau_v$ of all vertices decorated by 
$\bold{b}$ and respective orderings.

\smallskip

In the flowchart interpretation, such a corolla 
transforms an input family $(f_1,\dots,f_r)$,
$f_i:\,(\bold{Z}^+)^{a_{\bullet}} \to (\bold{Z}^+)^{c_i}$, into
the map 
$$
\langle f_1,\dots ,  f_r\rangle :\,
(\bold{Z}^+)^{a_{\bullet}} \to(\bold{Z}^+)^{c_1+\dots +c_r}.
$$
It was called juxtaposition in [Ma1], V.2.3 (b).
\medskip

(c) {\it Recursion.}  If a vertex is decorated
by $\bold{r}$, it must have exactly {\it two local inputs.}
If the arity/coarity of the first one (in their structure order) is 
$(a,c)$,
for the second one it must be $(a+c,c)$, and 
for the local output it must be $(a+1,c)$: this is our compatibility
condition.

\smallskip

In the flowchart interpretation, 
such a vertex takes as input two arbitrary 
maps $f_1:\,(\bold{Z}^+)^a\to (\bold{Z}^+)^c$,
$f_2:\,(\bold{Z}^+)^{a+c}\to (\bold{Z}^+)^c$
and produces the output 
$$
g:\,(\bold{Z}^+)^{a+1}\to (\bold{Z}^+)^c
$$
defined recursively as
$$
g(x,1):=f_1(x),
$$
$$
g(x, k+1):= f_2(x, g(x,k))
$$
for  each $x\in (\bold{Z}^+)^a,\, k\in \bold{Z}^+$.

\smallskip

This form of recursion is more restrictive than the one which is 
often used: it does not allow  $f_2$ to depend explicitly
on the recursion parameter $k$. However,
R.~M.~Robinson has proved in [Ro] that it suffices
to use it in order to get all primitive recursive functions,
if an extension of the list of basic functions is allowed.
Afterwards, M.~D.~ Gladstone  in [Gl] has shown that such an extension is unnecessary.
\medskip

{\it General case.} At first, consider a connected graph $\tau$.
Assume that it has $\ge 2$ vertices.
We define the operation
$Op\,(\tau )$ by induction on the number of vertices. 

\smallskip

Namely,  for a vertex $v$ which is the boundary of a global input, consider
the subfamily $f_v:= \{f_i\,|\,\partial_{\tau} (i)=v\}$.
Denoting by $\tau_v$ the corolla of $v$ (an {\it in--corolla}), calculate
$g_v:=Op\,(\tau_v) (f_v)$ as specified above.

\smallskip

One can check that this prescription produces the result independent
on arbitrary choices.

\smallskip

Now, consider the maximal decorated subgraph 
$\tau^0$ of $\tau$,
whose flags and vertices do not belong to this
$in$--corolla. Its global inputs consist of
all global inputs of $\tau$ not adjacent to $v$, and $j_{\tau}(r)$,
where $r$ is 
the root of our corolla.
Decoration of $\tau^0$ is the restriction of
that of $\tau$; global inputs of $\tau$ retain also
their input functions $f_i$. Decorate the input $j_{\tau}(r)$
by $g_v$ and put
$$
Op\,(\tau )(\{f_i\}):= Op\,(\tau^0)(\{f_i, g_v\,|\, \partial (i)\ne v\}).
$$
The right hand side is defined due to the inductive assumption.

\smallskip

Finally, if $\tau$ is the disjoint union of connected components
$\coprod_{a\in A}\tau_a$, we put
$$
Op\,(\coprod_{a\in A}\tau_a):= \times_{a\in A} Op\,(\tau_a)
$$
in the sense that $Op\,(\tau )$ acts on the family,
naturally  indexed by 
$A$, of (families of) global inputs of  connected components,
and produces the family of outputs, as well indexed naturally by $A$. 
  
\medskip

As we implied in the previous discussion, we can apply 
$Op\,(\tau)$ to families, consisting non--necessarily of basic,
or even recursive, functions.

\smallskip

But if we want to define programming methods
based upon $Prim$, then we must decorate global inputs by some basic functions, and interpret the resulting decorated tree as
as a program producing one concrete recursive function.

\smallskip

Here the choice becomes ambiguous:
we may change the list of basic functions, and we may allow
the application of $\bold{c},\bold{b},\bold{r}$
to some restricted class of subfamilies, getting the more general cases from trees larger than corollas.

\smallskip

For $\bold{c}$ and $\bold{b}$, we  allowed  arbitrary natural
families, implicitly using associativity
of intended interpretations. Yanofsky allows only
two inputs. For $\bold{r}$, 
we essentially adhered   to the choice made by
Yanofsky.

\medskip

{\bf  2.10. $Prim$ as a world of programming methods.} We may consider
$(\Z)^a$, $a\ge 0$, as objects of a category $\Cal{C}$ whose morphisms
are  recursive  functions. There are many other enumerable sets $S$ endowed
with a computable bijection (enumeration) $\Z \to S$. For example,
consider some finite Bourbaki structure, such as a finite group, or a finite graph.
One can easily construct a set of representatives
of isomorphism classes of such a structure, one from each class,
together with its enumeration.

\smallskip
Given such a set $\Cal{S}$,
we can add it to $\Cal{C}$ as a new object, with evident morphisms.
Such objects may be called (infinite) {\it constructive worlds.} For a more detailed
discussion, see [Ma2] and the 2nd Edition of [Ma1].
The new category will be equivalent to $\Cal{C}$. We may always assume
it to be monoidal with respect to the direct product $\times$.  

\smallskip

Now, the morphism sets $\Cal{C}((\Z)^a,(\Z)^c)$ that is,
recursive functions of a fixed arity and coarity, {\it do not form
a constructive world}: only their descriptions do.
Whenever a constructive world of descriptions is chosen, 
its definition must be completed by a family of {\it evaluation
maps} and {\it composition maps} which must be 
computable, that is, morphisms in the (extended) category $\Cal{C}$. 
When this is done, we will  elevate the constructive
world of descriptions to that of {\it  programming methods.} 

\smallskip

In order to illustrate these general considerations, define $P\,(a,c)$ as the subset of $Prim$, consisting of (isomorphism classes of)
graphs whose outputs (roots of connected components)
have the total arity/coarity $(a,c)$.

\smallskip

The evaluation morphism in $\Cal{C}$
$$
\roman{ev}_{P(a,c)}:\, P(a,c)\times (\bold{Z}^+)^a \to  (\bold{Z}^+)^c
$$
we have already essentially described. Namely,
$$
\roman{ev}_{P(m,n)} (\tau , (x_1,\dots ,x_m)):=
f_{\tau}(x_1,\dots ,x_m)
$$
where $f_{\tau}$ is the total  output of the flowchart $\tau$
which we formerly denoted $Op\,(\tau )$, applied to the input decorations
of $\tau$.

\smallskip

A computable multiple composition morphism 
$$
\roman{comp}:\ P(m_{r-1},m_r)\times \dots \times  P(m_2,m_3)\times P(m_1,m_2) \to
P(m_1,m_r)
$$
can be constructed as follows. For simplicity,
we will only describe the composite
$\roman{comp}\,(\tau_r,\tau_{r-1},\dots,\tau_1 )$
for an $r$--tuple of  decorated trees $\tau_1, \tau_2,\dots ,\tau_r$.  

\smallskip

 Consider a corolla with vertex decorated by $\bold{c}$,
$r$ inputs decorated by the arities $(m_1,m_2),\dots ,
(m_{r-1},m_r)$, and an output decorated by $(m_1,m_r)$.
Graft  inputs of ths corolla to the roots of $\tau_1,\dots ,\tau_r$ respectively. 
The resulting tree represents the composition. 

\smallskip

Of course, 
on the combinatorial level, we will have to make a stupid choice
of some ``concrete'' vertex and flags of this corolla,
 but the result will be unique up to unique
isomorphism identical on the component trees $\tau_i$.
\smallskip

However, if we iterate partial compositions which
on the level of maps correspond, say, to $h\circ g\circ f$,
$(h\circ g)\circ f$ and $h\circ (g\circ f)$ respectively, we will get
three different decorated trees, say $\sigma_{123}, \sigma_{12,3}, 
\sigma_{1,23}$. 

\smallskip

On the combinatorial/geometric level these trees are interconnected by two
contraction morphisms: $\sigma_{12,3}\to\sigma_{123}$
and  $\sigma_{1,23}\to\sigma_{123}$ which contract the edges
entering  to the root vertices, whose ends are marked by $\bold{c}$. 
One can simply declare that such contractions generate an equivalence relation on the elements of $Prim$, and that algorithms
encoded by $Prim$ are actually such (or even bigger) 
equivalence classes
rather than isomorphism classes of the decorated trees. 

\smallskip

However, since we work in a categorical context,
 a better way to proceed is to organize $Prim$ into
a constructive category, and then to {\it localize} it
with respect to those morphisms $\tau \to \sigma$ that  produce a natural identification $Op\,(\tau )$ and $Op\,(\sigma ).$ 

\smallskip
Recall
that the {\it localization of a category $\Cal{B}$ 
with respect to a set of its morphisms $S$}
is a functor $L:\,\Cal{B}\to \Cal{B}[S^{-1}]$ which makes all morphisms
in $S$ invertible and which is the initial object among all functors with this
property.  
 
\smallskip

Here is  a simple version of this construction.

\medskip

{\bf 2.11. Definition -- Claim.} {\it Consider the category $Pr$
whose set of objects is the set  $Prim$,
and morphisms are compositions of the following
maps of decorated graphs:

\smallskip

(i) Isomorphisms.

\smallskip

(ii) Contractions of subtrees of the following type:
all vertices of such a subtree are decorated by $\bold{c}$.
After the contraction, the resulting vertex must be marked by $\bold{c}$.
The remaining decorations do not change.
\smallskip

(iii) Contractions of subtrees, whose
all vertices are decorated by $\bold{b}$.
After the contraction, the resulting vertex must be marked by $\bold{b}$.
The remaining decorations do not change.

\smallskip

Denote by $P$ the  localization of $Pr$ with respect to all
morphisms. It has the natural structure of a category of programming methods for which composition and bracket operations become associative.}

\medskip

One can similarly accommodate more sophisticated equivalence
relations between decorated trees, studied by Yanofsky.

\smallskip

To this end one can extend  the category $Pr$ by some extra morphisms, 
and then localize with respect to them as well.  

\bigskip

\centerline{\bf 3. Bialgebras and Hopf algebras related to computation}
 
\medskip

Graphs which in the Feynman formalism enumerate terms
of a perturbation series, for the purposes of renormalization
are organized into another structure: they serve as formal
generators of an algebra endowed
with a diagonal map, that is, a bialgebra. The nature of this diagonal map $\Delta$
becomes even more transparent when graphs are treated as flowcharts: 
$\Delta$  simply
sums the different ways of decomposing  a flowchart into
a ``sub--'' and a ``quotient--''  component. Thus, a general context
for introducing the relevant bialgebras is that of morphisms in
categories, enriched categories, and programming methods.

\smallskip

In this section, we give a sample of such constructions.
The most basic and elementary example is that of finite categories.

 \medskip
 
{\bf 3.1.  Proposition.} {\it Let $C$ be a finite category.
Let $k$ be a base field (or unital commutative ring).
Denote by $B:=B_C$  the symmetric algebra of the linear space
freely generated by all morphisms in $C$. Define
the diagonal morphism of algebras $\Delta :\,B\to B\otimes B$
on generators:
$$
\Delta (f) :=\sum_{g,h\,|gh=f} h\otimes g
\eqno (3.1)
$$ 
Then $\Delta$ is coassociative so that $B$ becomes a bialgebra with
unit and counit.}

\medskip

{\bf Proof.} We have
$$
(\Delta\otimes \roman{id})\circ \Delta (f) = 
\sum_{k,l\,|kl=h}\left( \sum_{g,h\,|gh=f} l\otimes k\otimes g\right),
\eqno(3.2)
$$ 
$$
(\roman{id}\otimes\Delta )\circ \Delta (f) = 
\sum_{d,e\,|de=g}\left( \sum_{g,h\,|gh=f} h\otimes e\otimes d\right).
\eqno(3.3)
$$ 
It remains to check that each term at the right hand side of (3.2)
appears exactly once in the rhs of (3.3), and vice versa.
In fact, both sets are in a natural bijection
with triple decompostions $f=f_1f_2f_3$. This ends the proof.

\smallskip

To put it simply, coassociativity of $\Delta$ is a formal 
consequence of associativity of composition of morphisms.

\medskip

Composition of programming methods 
often is not strictly associative. Moreover, relevant  categories are not
always finite. Even Boolean circuits computing maps between finite sets
of a restricted cardinality do not form a finite set.

\smallskip

The latter problem can be alleviated by an additional structure that is always
relevant for programming methods: some natural $\Z$--measure
of {\it program size} and/or {\it computation time}, which is more or less
additive with respect to decompositions of programming methods
(as opposed to the decomposition of the respective functions) $F=GH$. 
Such an additivity  ensures that the right hand side of
formulas similar to (3.1) consists of a finite number of terms.
Associativity/coassociativity then is achieved by
appropriate modifications of constructive objects associated 
with programming methods.

\smallskip

The following example, based on graphs, displays some relevant constructions.
It is a version of Kreimer--Connes Hopf algebra in quantum field  theory
which can be re--interpreted as a composition bialgebra (3.1)
of programming methods as soon as an appropriate category of decorated
graphs is chosen.

\medskip

{\bf 3.2. Definition.} {\it Let $\tau$  be an oriented graph. Call {\it a proper cut} $C$
of $\tau$ any partition of $V_{\tau}$ into a disjoint union of two non--empty subsets 
$V_{\tau}^C$ (upper vertices) and $V_{\tau,C}$ (lower vertices) satisfying the following conditions:  

\smallskip

(i) For  each oriented wheel in $\tau$, all its vertices belong either to $V_{\tau}^C$,
or to $V_{\tau, C}$.

\smallskip

(ii) If an edge $e$ connects a vertex $v_1\in V_{\tau}^C$ to 
$v_2 \in V_{\tau,C}$, then it is oriented from $v_1$ to $v_2$
(``information flows down'').

\smallskip

(iii) There are also two improper cuts:  the upper improper cut
is the partition  $V_{\tau}^C=\emptyset $, $V_{\tau ,C}=V_{\tau}$,
whereas the lower one is
the partition  $V_{\tau}^C=V_{\tau} $, $V_{\tau ,C}=\emptyset .$}

\medskip

\smallskip 

Having chosen a cut $C$, we may define two graphs:
$\tau^C$ (upper part of $\tau$ wrt C) and  $\tau_C$ (lower part wrt C) by the following conditions. 

\medskip 
Vertices of   $\tau^C$  (resp. $\tau_C$) are $V_{\tau}^C$
(resp. $V_{\tau ,C}$). Flags of   $\tau^C$  (resp. $\tau_C$) are all flags
of $\tau$ incident a vertex of   $\tau^C$ (resp. $\tau_C$).
Edges of   $\tau^C$  (resp. $\tau_C$) are all edges
of $\tau$ whose both boundaries belong to   $\tau^C$ (resp. $\tau_C$).

\smallskip

Finally, all orientations remain the same as they were on $\tau$,
and if $\tau$ was not only oriented by additionally decorated,
or labels remain the same.
\smallskip

The term ``cut'' is motivated by the following simple observation. 
Given $(\tau ,C)$, we may define one more graph $\sigma$, with 
$(F_{\sigma}, V_{\sigma},\partial_{\sigma})$=
$(F_{\tau}, V_{\tau},\partial_{\tau}).$ Furthermore,
all edges remain the same except for those that lead from a vertex
in $\tau^C$ to a vertex in $\tau_C$: their halves become
some global outputs (resp. inputs) 
of $\tau^C$ (resp. $\tau_C$).

\smallskip

Thus, we have
$$
\sigma = \tau^C\coprod \tau_C
$$
where $\coprod$ is the disjoint union.

\smallskip

In other words, we get $|\sigma |$ by cutting all edges leading from
$|\tau^C|$ to  $|\tau_C|$ at their midpoints. This implies that $C$ can
be also identified with the respective subset of edges in $|\tau |$:
this is the definition of a cut often accepted in physics literature.

\medskip

{\bf 3.3.  Bialgebras of decorated graphs.}
Fix a set of labels $L=(L_F,L_V)$. Assume that $L_F=L^0_F\times \{in,out\}$,
When we speak about orientation of an $L$--decorated graph, we always 
mean the one that is obtained by forgetting labels from $L^0_F$.
The isomorphism class of a decorated graph $\tau$ is denoted $[\tau]$.

\medskip

{\bf 3.3.1. Definition.}  {\it A set $Fl$ (``flowcharts'')  of  $L$--decorated 
graphs is called  admissible,
if the following conditions are satisfied:

\smallskip

(i) Each connected component of a graph in $Fl$ belongs to $Fl$.
Each  disjoint union of a family of
graphs from $Fl$ belongs to $Fl$.  Empty graph $\emptyset$ is in $Fl$.

\smallskip

(ii) For each $\tau \in Fl$ and each   cut $C$ of $\tau$,
$\tau^C$ and $\tau_C$ belong to $Fl$.}

\medskip

Let now $Fl$ be an admissible set of graphs, and $k$ a commutative ring.
Denote by $H=H_{Fl}$ the $k$--linear span of isomorphism classes
of graphs in $Fl$: the $k$--module  of formal finite linear combinations $\{\,\sum_{\tau \in Fl} a_{[\tau]}[\tau]\,\}$.

\smallskip

Define two linear maps 
$$
m :\, H\otimes H\to H,\quad  \Delta :\, H\to H\otimes H
$$
by the following formulas extended by $k$--linearity:
$$ 
m ([\sigma ]\otimes [\tau ]):= [\sigma \coprod \tau ],
\eqno(3.4)
$$
$$
\Delta ([\tau ]) := \sum_C [\tau^C]\otimes [\tau_C],
\eqno(3.5)
$$
where the sum is taken over all cuts of $\tau .$

\medskip

{\bf 3.3.2. Proposition.} {\it (i) $m$ defines on $H$ the structure of a commutative 
$k$--algebra with unit $[\emptyset ]$. Denote the respective ring homomorphism 
$$
\eta :\,k\to H,\, 1_k\mapsto [\emptyset ]\, .
$$

\smallskip

(ii) $\Delta$ is a  coassociative comultiplication on $H$,
with counit 
$$
\varepsilon :\,H\to k,\ \sum_{\tau \in Fl} a_{[\tau ]}[\tau] \mapsto  
a_{[\emptyset ]}
\eqno(3.6)
$$

\smallskip

(iii) $(H,m,\Delta ,\varepsilon ,\eta )$ is a commutative bialgebra
with unit and counit.}

\medskip

{\bf Proof.} (i) The first statement is straightforward. It is worth mentioning
that $(H,m,\eta )$ is in fact the symmetric algebra freely generated
by the set of isomorphism classes $Fl_{con}$ of connected non--empty
graphs in $Fl$: $[\coprod_i \tau_i]$ corresponds to $\prod_i [\tau_i]$.

\smallskip

(ii) The least obvious in this statement is the coassociativity of $\Delta$. 
Omitting for brevity square brackets indicating isomorphism classes at the rhs,
we can write:
$$
(\Delta\otimes{id})\circ\Delta([\tau ])=\sum_{C}\sum_{C^{\prime}}
(\tau^C)^{C^{\prime}}\otimes (\tau^C)_{C^{\prime}}\otimes\tau_C,
\eqno(3.7)
$$
$$
(id \otimes\Delta)\circ\Delta([\tau ])=\sum_{C}\sum_{C^{\prime\prime}}
\tau^C\otimes (\tau_C)^{C^{\prime\prime}}\otimes(\tau_C)_{C^{\prime\prime}},
\eqno(3.8)
$$
where $C$ runs over cuts of $\tau$, $C^{\prime}$ runs over cuts of $\tau^C$,
and   $C^{\prime\prime}$ runs over cuts of $\tau_C$.

\smallskip

We want to establish a bijection between the sets of tensor monomials
in the rhs of both formulas.

\smallskip

To this end, consider triple partitions of $V_{\tau}$,
$$
V_{\tau} = V_1\coprod V_2\coprod V_3
$$
satisfying the conditions similar to those in Definition 3.2:

\smallskip

(a) For  each oriented wheel in $\tau$, all its vertices belong to one of the sets
$V_i$.

\smallskip

(b) If an edge $e$ connects a vertex $v_1\in V_i$ to 
$v_2 \in V_j$,  $i<j$, then it is oriented from $v_1$ to $v_2$.

\smallskip

From such a triple partition, we can produce two double partitions:
$(V_1\coprod V_2)\coprod V_3$ and  $V_1\coprod (V_2\coprod V_3).$
Both of them satisfy conditions of Definition 3.2. Hence they define two cuts of
$\tau$, say $C_{12}$ and $C_{23}$.

\smallskip

Moreover,  $V_1\coprod V_2$ defines a cut of $\tau^{C_{12}}$, say 
$C_{12}^{\prime}$, and hence a term in the rhs of (3.7). Similarly,
$V_2\coprod V_3$ defines a cut of  $\tau_{C_{23}}$, say 
$C_{23}^{\prime\prime}$, and hence  a term in the rhs of (3.8).

\smallskip

We claim that this construction establishes a bijection between the 
respective terms. The reasoning is somewhat cumbersome, but straightforward,
and we leave it to the reader.

\smallskip

The coidentity axiom reads 
$$
(\varepsilon\otimes{id})\circ\Delta([\tau ])=
(id \otimes\varepsilon)\circ\Delta([\tau ])= [\tau ].
\eqno(3.9)
$$
To check it, we refer to (3.5) and (3.6): only two terms
in (3.5), corresponding to improper cuts, can contribute to
(3.9). One gives the identity for the left counit, another for the
right one. 

\smallskip

(iii) It remains to check that $\Delta$ and $\varepsilon$
are algebra homomorphisms. For $\varepsilon$, this follows
from (3.4) and (3.5). For $\Delta$, this follows from the fact
that the a cut $C$ of $\sigma\coprod\tau$ is the same as
a pair of cuts  $(C_{\sigma},C_{\tau})$ of $\sigma ,\tau$ respectively,
so that $(\sigma\coprod\tau )^C= \sigma^{C_{\sigma}}\coprod
\tau^{C_{\tau}}$ and similarly for lower parts.

\medskip

{\bf  3.4. Hopf algebra of decorated graphs.} In order to construct
an antipode on the bialgebra $H=H_{Fl}$ we will show that 
one can introduce on $H_{Fl}$ a grading by $\bold{N}$ turning it into
a {\it connected graded bialgebra} in the sense of [E--FMan],
2.1. Then the existence of an antipode (and an explicit construction
of it) is provided by the Corollary 1 in 2.3 of [E--FMan].

\smallskip

There are two kinds of natural gradings. One can simply define
$$
H_n:=\ the\ k-submodule\ of\ H\ spanned\ by\ [\tau ]\ in\ Fl\  with\ |F_{\tau}|=n. 
$$
One can also introduce a weight function on the set of labels in 2.4:
$|\cdot|\,:\, L\to \bold{N}$  and put
$$
H_n:=\ the\ k-submodule\ of\ H
$$
$$
 spanned\ by\ [\tau ]\ in\ Fl\  with\ n=\sum_{f\in F_{\tau}}
(|l(f)|+1) + \sum_{v\in V_{\tau}} |l(v)|\, , 
$$
where $l:\, V_{\tau}\coprod F_{\tau} \to L$ is the structure decoration of $\tau$.
\smallskip

From the definitions, it follows that for either choice we have
$$
m(H_p\otimes H_q) \subset H_{p+q},\quad \Delta (H_n)\subset \oplus_{p+q=n} H_p\otimes H_q,
$$
and moreover,
$H_0 =k[\emptyset ]$ is one--dimensional, so that  $H$ is connected.

\medskip

{\bf 3.5. Hopf algebras from quantum computation: a review.} The last subsections of this part are dedicated
to one more class of constructions that are directly related to some ideas
in the quantum computation  project.

\smallskip

One standard model of quantum computation starts with a classical Boolean circuit $B$ which
computes a map $f:\,X\to X$, where $X$ is a finite set of Boolean words, say, of length $n$.
After replacing bits with qubits, and $X$ with the $2^n$--dimensional Hilbert space $H$
spanned by the ortho--basis of binary words, we have to
calculate the linear operator $U_f:\,H\to H$, linearization of $f$. Only unitary
operators $U_f$ can be physically implemented. Clearly, $U_f$ are unitary only
for bijective $f$; moreover,  they  must be calculated by ``reversible''  Boolean circuits.  

\smallskip
On the other hand, interesting $f$ are only rarely permutations.
For example, in search problems $f$ is the characteristic function of a subset
$X_0\subset X$ (``a needle  in a haystack'').

\smallskip

There is a well--known trick, allowing one to transform any Boolean circuit $B_f$
calculating $f$ into another Boolean  circuit $B_F$ of comparable length
consisting only of reversible gates and calculating a bijection $F$ of another finite set,
such that information about $f$ is easily read off from the corresponding 
information about $F$ (for more details, cf. [Ma2], 3.2).

\smallskip

If we now focus on permutations of $X$, there naturally arise two Hopf algebras
related to them: {\it group algebra of permutations} and a dual Hopf algebra.
For infinite $X$, there are several versions of Hopf algebras
associated to the group of unitary operators $H_X\to H_X$.

\smallskip

Below we reproduce the combinatorial part of these constructions,
having in mind applications to quantum computations of
{\it recursive} functions $f$, say, $\Z\to\Z$.
An additional complication arises here: $f$ can be only partial,
hence before processing  it into a permutation, we must make
it everywhere defined. 

\smallskip

Before proceeding to details, notice a change of perspective.
Working with $Prim$ as flowcharts in 2.9, and by extension in 3.3,
 we considered {\it programming
methods}, whose inputs and outputs were functions.
Here we work with what earlier was called {\it instances of computation:}
inputs and outputs are now arguments/values of a function to be computed.

\medskip

{\bf 3.6.  Reduction of total maps to bijections.} Consider a 
set $X$ and a class $\Cal{F}$ of 
everywhere defined maps $f:\, X\to X$.

\smallskip

In order to reduce $\Cal{F}$ to permutations,
introduce on $X$ a structure of, say, abelian group with composition law 
denoted $+$. Produce from $f$ the map
$$
\tilde{f}:\, X^2\to X^2, \ \tilde{f}(x,y):=(x+f(y),y).
\eqno(3.10)
$$ 
This is a bijection: $\tilde{f}^{-1}$
maps $(x^{\prime}, y^{\prime})$ to $(x^{\prime}-f(y^{\prime}),y^{\prime})$.
Knowing $\tilde{f}$, we can compute $f$: take the first coordinate
of $\tilde{f}(0,y)$. 

\smallskip

If $X$ is endowed with a natural enumeration, one should choose an (easily) computable group law $+:\,X^2\to X$, and thus reduce the computation of $f$
to that of $\tilde{f}$ using a (hopefully manageable) additional amount of memory and time.  In turn, if $f$ is computable, $\tilde{f}$ will be as well.
 Identical permutations  and composition of two computable permutations are
 computable.

 \medskip

{\bf 3.7. Reduction of partial maps to total maps.} First recall how one composes partial
maps.
\smallskip
Formally, a partial map from a set $X$ to a set $Y$
is  a pair $(\varphi,D(\varphi))$ where $D(\varphi)$ is a subset of $X$ (possibly  empty),
and  $\varphi:\,D(\varphi )\to Y$ is an actual map. We put $\roman{Im}\,\varphi :=
\varphi(D(\varphi )).$
Denote $Par\, (X,Y)$ the set of partial maps.
The composition $Par\, (Y,Z)\times  Par\, (X,Y)\to Par\, (X,Z)$ is defined  as 
$$
(\chi ,D(\chi ))\circ (\varphi ,D(\varphi )) := (\chi \circ\varphi,\varphi^{-1}(D(\chi)\cap \roman{Im}\,\varphi)\,).
$$
One easily sees that in this way we get a category, say $ParSets.$

\smallskip

Notice that each set of morphisms $Par\,(X,Y)$ is  {\it pointed},
in the sense that it has a canonical element, ``empty map'', say,
$\emptyset_{X,Y}$. Its composition with any other morphism is again
the respective empty map.

\smallskip

 This last remark motivates the consideration of another
category: that of {\it pointed sets}  $PSets$. An object of
$PSets$ is a pair $(X,*_X)$ where $*_X\in X$ (so that $X$
cannot be empty). A morphism $(X,*_X)\to (Y,*_Y)$ is an everywhere defined
map $f :\, X\to Y$ such that $f (*_X)=*_Y.$
The composition is evident.

\smallskip

Deleting marked points, we get a functor  $PSets\to ParSets$:
$$
X\mapsto X^{\circ}:= X\setminus \{*_X\},\quad f \mapsto f^{\circ}:=(\varphi, D(\varphi)),
\eqno(3.11)
$$
where, for $f:\,X\to Y$, $D(\varphi )$ is defined as  $f^{-1}(Y^{\circ})$
and $\varphi$ as the restriction of $f$ to $D(\varphi )$.

\smallskip

This construction is obviously invertible in the sense that
there exists  a quasi--inverse functor $ParSets\to PSets$.
It can be constructed by formally
adding an extra marked point  $*_X$ to each object $X$
in $ParSets$, and extending each partial map $(\varphi ,D(\varphi ))$
from $X$ to $Y$ by sending $X\setminus D(\varphi )$ to $*_Y$.

\smallskip

Let us return now to the situation where $X$, $Y$ are endowed 
with computable numberings, and restrict ourselves to 
semi--computable partial functions $\varphi$. Then the passage from $\varphi$
to an everywhere defined (total) function $f:\,X\cup\{*_X\}\to Y\cup \{*_Y\}$ involves 
a simple extension of $\varphi$: we put  
$$f(x)=\varphi (x)\  \roman{for\ all}\ x \in D(\varphi ),\
f (x)=*_Y\ \roman{otherwise.}
\eqno(3.12) 
$$

\smallskip

Of course, $Y\cup\{*_Y\}$ is endowed with obvious computable
numberings  compatible with that of $Y$, say, one can simply augment by 1 
the initial numbering of $Y$ and put $*_Y$ at the first place. But from the
viewpoint of computability  $*_ Y$ looks rather as ``infinite'' , or ``transfinite''
element: if $D(\varphi )$ is only enumerable but not decidable, a Turing machine trying to calculate  
$f(x)$ for $x\in D(\varphi )$ might never stop. For the same reason,
$f$ as a total function might become uncomputable.

\smallskip

Nevertheless, we can apply to $f$ the trick (3.10) and get a
 permutation $\tilde{f}$ of  $(X\cup\{*_X\})^2$ which will be uncomputable
 outside $(X\cup\{*_X\})\times (D(\varphi )\cup \{*_X\})$. Choosing a computable
 structure of abelian group on $X\cup\{*_X\}$ we will have to treat $*_X$
 as an ordinary element. Choosing it to be zero, we
 get the following nice statement.
 
 \medskip

 {\bf 3.8. Proposition.} {\it Let $\varphi :X\to X$ be a partial recursive function.
Construct its extension as above $f:\,X\cup\{*_X\}\to X\cup \{*_X\}$.
Choose a computable (general recursive)  structure of additive group on 
$X\cup \{*_X\}$ with zero $*_X$.

\smallskip
 
 Denote by $\tilde{f}$ the permutation of $Z:=(X\cup\{*_X\})^2$
 produced from $f$  by  the analog of formula (3.10) in this context. Then 
 $\tilde{f}$ is a permutation with the following
 properties:
 
 \smallskip 

(i) $\tilde{f}$ is an extension of the partial recursive function $g :\,Z\to Z$
with the definition domain
$$
D(g ):= (X\cup\{*_X\})\times (D(\varphi )\cup\{*_X\}).
\eqno(3.13)
$$

(ii) $\tilde{f}$ induces a permutation of $D(g )$ with unique fixed point
$(*_X,*_X)$. The
 complement of $D(g )$ consists  of all remaining fixed points of $\tilde{f}$.
}

\medskip

{\bf Proof.} By definition, we have
$$
x\in D(\varphi )\  \Longrightarrow  \  f(x)=\varphi (x),
$$
$$
x\in X\cup\{*_X\}\setminus D(\varphi ) \Longrightarrow  \  f(x)= *_X.
$$ 
Therefore,  denoting by $x,y$ variable elements of $X\cup\{*_X\}$, we have:
$$
\tilde{f}(x,y) = (x+\varphi (y),y),\ \roman{if} \ y\in D(\varphi ),
$$
$$ 
\tilde{f}(x,y) = (x,y),\ \roman{if} \ y\notin D(\varphi ),
$$
Therefore, $\tilde{f}$ is computable on (3.13), and has there a unique
fixed point $(*_X,*_Y)$. 
Moreover, since  in the  case $y\in D(\varphi )$,  $\varphi (y)$ is never zero, we easily obtain (ii). 

\bigskip

\centerline{\bf 4. Regularization and renormalization}

\medskip

The subsections 4.1--4.5 are dedicated to a review 
of the relevant parts of the renormalization formalism
in Quantum Field Theory, following
[E--FMan]. The reader should keep in mind that this scheme does not cover
all versions, used by physicists: see e.g. [Cos] for a detailed treatment
of the so called Wilsonian renormalization, and [Po] for a brief introduction to the peculiarities
related to the quantization of gauge fields.

\smallskip

In the subsections 4.6--4.9 I review some regularization schemes that might be relevant
in theoretical computation.

\medskip

{\bf 4.1. Connected filtered Hopf algebras.} Let $\H$ be
a unital associative and counital coassociative bialgebra over a field $K$ 
of characteristic zero. The relevant structure maps are denoted
$$
m:\,\H\otimes\H\to \H,\ \Delta :\,\H\to \H\otimes \H ,\
u:\, K\to \H,\ \varepsilon :\, K\to \H .
$$ 
In our main applications $\H$ will be a bialgebra of programming methods,
such as flowcharts (see 3.3).
 It will satisfy two main assumptions
of [E--FMan]: it will be filtered and connected. This means that we are given
a filtration $\H=\cup_{n=0}^{\infty}\H^n$,  compatible with $m$ and $\D$
in the standard sense:
$$
m({\H}^p \otimes \H^q)  \subset \H^{p+q}, \quad \D (\H^n)\subset \sum_{p+q=n}\H^p\otimes\H^q,
$$
and  moreover, $\H^0$ is identified with $K$ by means of
$u$ and $\varepsilon$. 

\smallskip

In this case, $\H$ automatically has an antipode $S$ with $S(\H^n) \subset\H^n$, and hence is a Hopf 
algebra. The antipode can be given explicitly, by induction on the filtration degree.
Namely, we have $S(1)=1$, and for any $x\in \H^n,\,n\ge 1$, in a version of Sweedler's notation
$$
\widetilde{\D}x:= \D x-(x\otimes 1+ 1\otimes x)= \sum_{(x)}x^{\prime}\otimes x^{\prime\prime} \in \ \bigoplus {\Sb p,q\ge 1 \\ p+q=n \endSb} \H^p\otimes\H^q
\eqno(4.1)
$$
and
$$
S(x)=-x- \sum_{(x)}S(x^{\prime}) x^{\prime\prime} = -x
- \sum_{(x)}x^{\prime} S(x^{\prime\prime})
 \eqno(4.2)
$$

The antipode is a crucial ingredient in the renormalization formulas 
(4.6) and (4.7) below.

\smallskip

Since its existence is guaranteed by (4.1), (4.2),  in applications
to programming methods, when the relevant bialgebra is constructed,
we must define {\it a compatible filtration.} Often it comes from a
bialgebra  {\it grading,}  whose intuitive meaning is
quite transparent: it is a measure of complexity/volume of the
relevant programming method which is additive with respect
to the composition of programs.

\smallskip

For example, the total number of flags is additive with respect to cuts
(cf. Definition 3.2), so it can be used to define the grading
of the respective bialgebra. When graphs, such as flowcharts, are decorated by labels
from a countable set $L$ , one can choose a ``weight'' numbering of $L$
and define the grading degree of  a decorated  graph as {\it the sum of all
labels of its flags and vertices.} 

%In our applications, the following additional structure might be relevant.

%\medskip

%{\bf 1.1.1. Definition.} {\it A positive semilattice $\N \subset \H$ is a subsemiring
%of the form $N[x_a\,|\,a\in A]$ (polynomials of a family of free commutative
%generators with non--negative integral coefficients) such that
%$x_a\in \r{ker}\,\varepsilon$, and for any $x=x_a\in \N$, the elements
%$x^{\prime}, x^{\prime\prime}$ from (1.1) can also be taken in $\N$.}

%\smallskip

%For example, in QFT renormalization schemes, $x_a$ are isomorphism 
%classes of Feynman diagrams.

%\smallskip

%We will denote by $\N\otimes \N$ the polynomial  semiring
%$$
%\N\otimes \N:=N[\,x_a\otimes 1, 1\otimes x_a\,|\,a\in A\,].
%$$

\medskip
{\bf 4.2. ``Minimal subtraction'' algebras.}  A "minimal subtraction" scheme,
that merges well with Hopf algebra renormalization techniques (cf. [E-FMan]),
formally is based upon a commutative associative algebra $\A$ over a field $K$,
together with two subalgebras $\A_-$ ( ``polar part'')  and  $\A_+$ (regular part),
such that $\A=\A_-\oplus \A_{+}$ as a linear space. One usually assumes
$\A$ unital, and $1_A\in \A_+$.  Moreover, an augmentation homomorphism
$\varepsilon_{\A} :\,\A_+\to K$ must be given.

\smallskip

Then any element $a\in \A$ is the sum
of its polar part $a_-$ and regular part $a_+$. The ``regularized value'' of $a$
is  $\varepsilon_{\A} (a_+)$.  

\medskip

{\bf 4.3. Example: germs of meromorphic functions.} Here $K=\bold{C}$,
$\A :=$ the ring of germs of meromorphic functions of $z$, say, at $z=0$;
$\A_- :=z^{-1} \bold{C} [z^{-1}]$, $\A_+$ consists of germs of regular functions at $z=0$, 
$\varepsilon_{\A} (f):= f(0).$

\smallskip

Notice that for the same algebra,
a complementary choice could have been made: 
one could  put $\A^{\prime}_-:= $ germs regular and vanishing at $z=0$,
and $\A^{\prime}_+:= \bold{C}[z^{-1}]$, with $\varepsilon_{A}^{\prime}(f)=f(z_0 )$
for a constant $z_0$.

\smallskip

This is a toy model of  situations arising in  cut--off regularization schemes:
$z_0$ is an input  from physics,  the scale of a parameter $z_0$ 
(such as energy) to which our observed quantities refer. 
As a physics justification of the choice of such a scale one
might postulate the belief that beyond this scale
`` a new physics'' starts.

\medskip

{\bf 4.4. Hopf renormalization scheme.} We summarize its algebraic
version here, restating Theorem 1 from sec. 2.5, [E--FMan].

\smallskip

Let $\H$ be a Hopf algebra as above,  $\A_+,\A_-\subset \A$
a minimal subtraction unital algebra. Consider the set $G(\A )$ of $K$--linear maps
$\varphi :\,\H\to\A$ such that $\varphi (1_{\H})=1_{\A}$.

\smallskip

Then $G(\A )$ with the convolution product
$$
\varphi*\psi (x):= m_{\A}(\varphi\otimes \psi )\D (x)= 
\varphi(x) +
\psi (x)+ \sum_{(x)}\varphi (x^{\prime})\psi (x^{\prime\prime})
\eqno(4.3)
$$
is a group, with identity $e(x):= u_{\A}\circ \varepsilon (x)$ and inversion
$$
\varphi^{*-1}(x)= e(x)+\sum_{m=1}^{\infty} (e-\varphi )^{*m}(x)
\eqno(4.4)
$$
where for any $x\in\r{ker}\,\varepsilon$  the latter sum contains only finitely many
non--zero summands.

\smallskip

An important subset of $G(\A )$ consists of {\it characters:} those linear
maps $\Cal{H}\to \A$ that are homomorphisms of algebras.

\medskip

{\bf 4.4.1. Theorem on the Birkhoff decomposition.} {\it If $\A$
is a minimal subtraction algebra,  each $\varphi\in G(\A )$
admits a unique decomposition of the form
$$
\varphi =\varphi_-^{*-1}*\varphi_+;\quad \varphi_- (1)=1_{\A},\ \varphi_{-} (\r{ker}\,\varepsilon )
\subset \A_-,\ \varphi_+(\H)\subset \A_+.
\eqno(4.5)
$$
Values of renormalized polar (resp. regular) parts $\varphi_{-}$
(resp. $\varphi_+$) on $\r{ker}\,\varepsilon$ are given by the inductive formulas
$$
\varphi_-(x)=-\pi \left(\varphi(x) +
 \sum_{(x)}\varphi_- (x^{\prime})\varphi (x^{\prime\prime})\right),
\eqno(4.6)
$$
$$
\varphi_+(x)=(\r{id}-\pi ) \left(\varphi(x) +
 \sum_{(x)}\varphi_- (x^{\prime})\varphi (x^{\prime\prime})\right).
\eqno(4.7)
$$
Here $\pi :\,\A\to \A_-$ is the polar part projection in
the algebra $\A$.

\smallskip

If $\varphi$ is a character,   $\varphi_+$ and  $\varphi_-$ are characters as well.}

\medskip

{\bf 4.5. Rota--Baxter operators as generalized polar parts.}
The general definition
of a Rota--Baxter (RB) operator of weight $\theta$ on an  associative (not necessarily unital or commutative) algebra $\A$ is this: it is a linear operator $R:\,\A\to\A$,
satisfying the identity
$$
R(f)R(g)= R(R(f)g+fR(g)+\theta fg)
\eqno(4.8)
$$

If $\A$ is a minimal subtraction algebra, the polar part projection
$\pi :\,\A\to\A_-$ is an RB operator of weight $-1$; moreover,
$\pi^2=\pi$.  

\smallskip

Dropping the restrictions $\theta =-1$ and $R^2= R$,
but still imagining $R$ as a generalized  ``polar part''
operator, one gets more freedom in using the recursive renormalization
schemes (4.6), (4.7): see [E--FMan] and below.

\medskip

We now pass to the discussion of this scheme in possible
applications to programming methods. We have already explained that to construct 
the relevant Hopf algebra $\Cal{H}$ we need a composition bialgebra as 
in section 3, completed by a filtration coming, say, from an additive
complexity  function on programming methods.  

\smallskip

{\bf 4.6. Target algebras and tropical geometry.} 
We now turn to possible interpretations of  target algebras
$\Cal{A}$ and linear functionals/characters $\varphi :\,\Cal{H}\to \Cal{A}$ in the theory
of computation. 

\smallskip

Roughly speaking, there are at least two classes of meaningful pairs $(\A ,\varphi )$:
 
\smallskip

(A) $\varphi (x)$ can be a quantitative characteristic (measure) {\it of a  program $x$}
(e. g.  running time and/or memory needed to complete the computation as a function
on the size of input);

\smallskip
(B) $\varphi(x)$ is  quantitative characteristics {\it of the output}
produced by a variable program $x$ 
(on a particular input, or on the set of all inputs). 

\smallskip

Some natural measures take their values not in commutative rings,
but in {\it commutative semirings}, say, of the type {\it Max--Plus}: cf. [Cas1], [Cas2]
for a recent review and references. Such a structure is similar to a commutative ring,
but its additive group axioms are weakened to those of additive monoid.
The name {\it Max--Plus} is explained by a typical example of such a semiring:  
$\bold{R}_{\ge 0}$ with ``addition'' $\r{max}\,(x,y)$ and ``multiplication'' $x+y$.

\smallskip

Such measures do not fit directly in the framework of Hopf renormalization
theory as it was formulated, but I want to stress their importance here,
by providing two examples. I plan to return to Hopf semiring renormalization
in a sequel to this article. 

\medskip

{\bf 4.6.1. Example: parallel computation.}  Let a world $W$ of histories
of computations be represented by  decorated flowcharts  as in 3.3. 
 Then a reasonable idealization of
{\it running time} might be a function $T:\, W\to {\R}_{\ge 0}$ with the following
property: for any flowchart $\tau$ and its cut $C$,
$$
T(\tau )=T(\tau^C)+T^C(\tau_C)
\eqno(4.9)
$$
where the superscript in $T^C$ is supposed to remind that 
input of $\tau_C$ is the output of $\tau^C$.
On the other hand, the idea of {\it parallel computation} is reflected in a formal
requirement:  running time of a disjoint sum of flowcharts is 
$$
T(\tau_1 \coprod\tau_2) =\r{max}\,(T(\tau_1),T(\tau_2)).
\eqno(4.10)
$$
Such a function can be thought of as a semiring--valued  ``quasi--character'' 
sending the composition of programs to the semiring product, and 
disjoint sum of programs to the semiring sum.

\smallskip

For more details, see sec. 2 of [Ma4].
\medskip

{\bf 4.6.2. Example: Bayesian networks.} A Bayesian network (see e.~g. [PaSt1--2])
 is a directed graph
whose vertices are (decorated by) two groups of variables:
observable random variables $Y_1,\dots ,Y_n$ and hidden
random variables $X_1,\dots, X_m$, whereas edges are decorated by (matrices of)
transition probabilities. A Bayesian network can be considered as a programming method
for computing certain characteristics of hidden random variables, such as so called
{\it marginal probabilities}, and {\it maximal a posteriori $\r{log}$ probabilities (MAP)}
(cf. [PaSt1]).  Although Bayesian networks are not, strictly speaking, flowcharts in
the sense of  our sec. 3, it is not difficult to recast them in a similar form.
Then the outputs such as marginal and MAP probabilities become
semiring--valued characters. 

\medskip

{\bf 4.7. Algebras of sequences.}  We start in a formal setting.
For a field $K$, consider the space  $\Cal{S}=\Cal{S}_K$
of infinite sequences $f=(f_1,f_2,\dots ),\, f_i\in K.$ There are three
different multiplications that furnish relevant structures of commutative algebra on this space, 
$\b , \times$
and $*$:
$$
(f\b g)_n:= f_ng_n,\ 1_{\b} =(1,1,1, \dots ), 
\eqno(4.11)
$$
$$
(f* g)_n:= \sum_{\r{max} (p,q)=n} f_pf_q, \ 1_{*} = (1,0,0,\dots ),
\eqno(4.12)
$$
$$
(f\times g)_n:= \sum_{p+q=n} f_pf_q,\quad  \r{(non-unital).}
\eqno(4.13)
$$
\smallskip
Clearly, $(\Cal{S}_K,\times )$ is simply the algebra of formal series
$\sum_{n=1}^{\infty} f_nz^n$. We could  formally adjoin an identity to it,
that is, allow sequences starting with $n=0$.  Multiplications 
$\b$ and $*$ extend to this case without problems.
However, we will have to avoid non--vanishing $f_0$ in other contexts,
cf. below. 
\smallskip

We can similarly interpret $(\Cal{S}_K, *)$ as the algebra of formal
``tropical''  series
$\sum_{n=1}^{\infty} f_nz^n$, with multiplication $z^p*z^q:=z^{\r{max} (p,q)}$.
Notice that $z$ is the identity in this algebra. However, if we extend it
by allowing $f_0\ne 0$, the role of identity will pass  to the constant
formal series $1$.

\smallskip

The map ``partial summation'':
$$
S:\,(\Cal{S}_K,*) \to   (\Cal{S}_K,\b ),\quad S(f)_N:=  \sum_{n=1}^N f_n
\eqno(4.14)
$$
is an isomorphism of unital algebras:
$$
S(f*g)=S(f)\b S(g).
$$

\smallskip

Considered as a map of the algebra  $(\Cal{S}_K,*)$ {\it into itself}, $S$ is 
{\it a Rota--Baxter (RB) operator of weight 1} in the sense
of [E--FMan], Sec. 3: we have
$$
S(f)*S(g)= S(S(f)*g+f*S(g) +f*g).
\eqno(4.15)
$$
The role of this remark in the context of renormalization
is explained by the fact that in the simplest  ``minimal subtraction'' scheme as in
the subsection 1 above,
the projection to the ``polar part'' $\pi :\, \A\to \A_-$ is an idempotent 
Rota--Baxter operator of weight $-1$. 

\smallskip

Notice in conclusion that the slightly modified summation operator $S^{\prime}$,
$$
S^{\prime}(f)_N:= \sum_{n=1}^{N+1} f_n,
\eqno(4.16)
$$
is again Yang--Baxter, of the same weight $-1$ as a polar projection $\pi$ is.

\medskip
{\bf 4.7. Boutet de Monvel's regularization.} Let now $K=\C$.
Consider  the subspace
$\A =\A_{\bold{C}}$ in $\Cal{S}_{\bold{C}}$ consisting of such sequences  
$f=(f_n)$ for which there 
exists a polynomial $P=P_f$ and an integer $A =A_f$ with the property
$$
S(f)_N= P_f(\roman{log}\,N) + O((\r{log}\,N)^A/N)
\eqno(4.17)
$$   
\medskip

{\bf 4.7.1. Theorem.} {\it   (i) $S(\A )$ is a subalgebra in  
$(\Cal{S}_{\bold{C}},  \b )$, so that $\A$ is a subalgebra of tropic power series 
$(\Cal{S}_{\bold{C}},  * )$.

\smallskip

(ii) Considered as  a subset of  $(\Cal{S}_{\bold{C}},  \times) =z\bold{C}[[z]]$, $\A$ can be described as 
a subspace of the algebra $\Cal{B}$
of power series $f(z)$ convergent in $|z|<1$ and vanishing at 0, such that
for a polynomial $Q=Q_f$ and an integer $B=B_f$ we have as $z\to 1-0$
along the real axis
$$
f(z)= Q_f(-\r{log}\,(1-z))+O((1-z)\cdot \r{log}^B(1-z)).
\eqno(4.18)
$$
\smallskip

(iii) The ``singular part"  $Q_f$ is uniquely defined
by $P_f$. It  can be derived from $P_f$ by applying  the formal differential 
operator of infinite order obtained from the Taylor series of the gamma--function:
$$
Q_f(t)= \Gamma (1+\p_t) \,P_f(t),\quad \Gamma (1+x) = 1-\gamma x+\dots .
\eqno(4.19)
$$
}

\medskip
The first part is straightforward. The second and third ones constitute a theorem by
Boutet de Monvel: see [BdeM1] for a more precise statement, and [De],
[Ra] for some proofs.

\medskip

{\bf NB}  The embedding $\A\subset \Cal{B}$ is strict.

\bigskip

{\bf 4.7.2. Comments.} Consider the subalgebra $\Cal{B}_-\subset \Cal{B}$
consisting of all functions $Q (-\r{log}\,(1-z))$, $Q\in t\,\bold{C}[t]$.
In fact, we have $\Cal{B_-}\subset \A$, so we
may alternatively call it $\A_-$.

\smallskip
The map $\pi :\,f\mapsto Q_f$  is a surjective algebra homomorphism
$\Cal{B}\to \Cal{B}_- \subset \Cal{B}.$ 
It is a natural singular part of $f$ at $z=1$ replacing the polar
part  of a meromorphic function.
Its kernel $\Cal{B}_+$ (resp $\A_+$) consists of elements of $\Cal{B}$
(resp. $\A$) vanishing at $z=1$. 
Hence we get a minimal subtraction algebra $\bold{C}\oplus \Cal{B}$ that 
can  serve as an input
to a natural  regularization scheme.  In fact, Boutet de Monvel, Racinet
and Deligne used it to regularize
the integral expressions for multiple zeta values which are closely related
to special Feynman integrals.

\smallskip
Restricted to $\Cal{A}$, $\pi$ is a linear
surjective map $\Cal{A}\to \Cal{A}_-$.

\smallskip

As at the end of 4.3, we can also interchange the nominations of
``polar'' and ``regular'' parts.

\smallskip 

Finally,  notice that $-\r{log}\,(1-z)$ corresponds to the sequence
$$
\bold{l}:= (1, 1/2, 1/3, \dots ).
\eqno(4.20)
$$
This sequence satisfies  (4.17) with $P_f(t)=t+\gamma$.
Hence it is natural to declare that the subalgebra of polar parts
in $(\A, *)$ (i.~e. in $\A$ interpreted as tropical series)
is the algebra of $*$--polynomials in $\bold{l}$. The formula
(4.19) says that  the change of multiplication from $*$ to $\times$ does not
change the space of polar parts.

\medskip

Results of A.~Levin ([Le])  suggest  that in (one of the) computation contexts
one can meaningfully  replace (4.20) by a sequence related to
exponential Kolmogorov complexity, in the spirit of the remark made in 0.6 above
that Kolmogorov complexity is the ultimate computational infinity.
I will finish this paper by briefly explaining Levin's theorem.

\bigskip

{\bf 4.8. One--sided enumerability.} The standard definition of a computable
real number $x$ involves a recursive sequence of rational approximations $(r_n)$
to it, together with a recursive sequence of bounds for the error $|x-r_n| \le b_n$,
$b_n\to 0.$

\smallskip  

We will consider here, following [Le],   {\it one--sided versions}
of computability
that are related to various versions of Kolmogorov complexity.

\medskip

{\bf 4.8.1. Definition.} {\it A real number
$x\in \bold{R}$ is called  enumerable from below iff  the following
equivalent conditions (i), (ii) are satisfied:

\smallskip
(i) there exists
a general recursive function $\varphi :\, \Z \to \Q$ such that
$$
\varphi (1)\le \varphi (2) \le \varphi (3) \le \dots,\quad \r{lim}\,\varphi (n)=r.
\eqno(4.21)
$$

\smallskip

(ii) The set $\{r\in \Q\,|\, r< x\}$ is recursively enumerable.}

\medskip

For obvious reasons, we consider the symbol $+\infty$ as computable from below as well
(this symbol is useful in  tropical contexts): it is approximated by $\varphi (n)=n$.

\medskip

In [BrYa] such numbers are called left computable (or rather, $-r$ are called right computable).

\bigskip

{\bf 4.8.2. Definition.} {\it A sequence of real numbers
$x_n$, $n=1,2,\dots $ is called enumerable from below iff  the following
equivalent conditions (i), (ii) are satisfied:

\smallskip
(i) there exists
a general recursive function $\varphi :\, \Z \to \Z\times \Q$, $\varphi(n)=(m_n, r_n)$ such that
the map $\Z\to \Z :\,n\mapsto m_n$ is surjective; if $m_a=m_b,\ a<b$, then
$r_a\le r_b$; and finally, the limit of the sequence $r_a^{(n)}$ corresponding to one and the same
first coordinate $n$, is $x_n$.

\smallskip

(ii) The set $\{(n,r_n)\in \Z\times \Q\,|\, r_n< x_n\}$ is recursively enumerable.}

\medskip

Again, we may include the symbol $+\infty$ as a possible value of $x_n$.

\medskip

We will use results of [Le], where such sequences of non--negative real numbers 
were considered.

\bigskip

{\bf  Remarks.}  (a) Enumerable from below reals form an additive subsemigroup of
$\R$. The check is straightforward. The same is true for reals enumerable from above (one can iclude the symbol 
$-\infty$ in place for $+\infty$).

\smallskip

(b) There exist reals enumerable from below but not  from above,
and vice versa. In fact, if a real $x$ is enumerable both from below and above,
than it is computable: the differences between the upper and the lower
$n$--th approximation form a recursive sequence of bounds for each of
the approximations.

\smallskip
It is crucial that  there exist one--sided  enumerable numbers which are not computable:
see e.~g. Proposition 2.2 in [BrYa].  

\smallskip

(c) For this reason, enumerable from below reals {\it do not form a subring
of $\R$}: multiplication, say, by $-1$  reverses enumerability from below
to enumerability from above. 

\smallskip

Similarly, {\it inversion} $x\mapsto x^{-1}$   interchanges
enumerability from below
and enumerability from above. 

\smallskip
However, {\it non--negative  reals} enumerable
from below (resp. from above) form a semiring with respect to
the usual addition and multiplication.

\smallskip

They also form a tropical semiring with respect to the operations 
$\oplus :=\r{max}, \otimes :=\cdot$ (resp. $\oplus :=\r{max}, \otimes :=\cdot$).

\medskip

{\bf 4.9.  Kolmogorov complexity reappears. } The simplest construction due to L.~Levin
starts with a particular  ``norm functional'' defined
on the set of enumerable from below sequences $f:\, \Z\to \bold{R}_{\ge 0}$
(our former $x_n$ is now $f(n)$):
$$
 N(f):=\r{sup}\,\{r\cdot \r{card}\,\{x\,|\,f(x)\ge r\}\}
$$ 

\medskip

{\bf 4.9.1. Proposition.}  {\it  There exists an enumerable from below sequence of finite norm $F$
such that any enumerable  from below sequence $f$ is majorized by $cF$
for an appropriate constant $c$.

\smallskip

The sequence $(-\r{log}\, F(n))$ coincides with the sequence of values of
(logarithmic) Kolmogorov complexity, up to an $O(1)$ function.}

\bigskip
l intend to study in the sequel of this paper a version of regularization where (4.20) is replaced by
$(F(n))$.

\bigskip

\centerline{\bf Appendix: Renormalization at large}

\medskip

{\it Warning: the gentle reader is kindly invited to skip the following musings.}

\medskip

{\bf A.1. Is physical infinity real?} The basic  intuitive picture
behind mathematical formalism of renormalization seems to be  
an image of   ``finite''  {\it observable} reality as 
{\it a difference (or quotient)} of two unobservable and infinite {\it physical} realities. The same intuition in theoretical physics produces such
expressions as ``vacuum energy is infinite'': the observable
finite lumps of energy/matter are interpreted as finite differences
between the the two excitation levels of vacuum, both of which
are infinite. 

\smallskip

This intuition  is supported by technological achievements:
the energy of a nuclear explosion is freed, when the two
infinities, tightly  balanced at the nuclear  scale, are made unbalanced by subtly
controlled technological processes -- their difference then destroys everything around them.

\smallskip

Finally, the question in the title of this subsection should not be
confused with a totally different question:  {\it ``Is physical reality infinite?''}
(Wittgenstein's angry laugh thunders from the 
Great Beyond  ...)

\medskip

{\bf A2. Epistemology of mathematics and infinite.} A possible parallel in the platonic world of ideas can be 
traced in the on--going epistemological shift
related to the foundations of mathematics: {\it discrete and finite}
 nowadays often comes
from looking at (homotopy types of) {\it continuous and infinite.}

\smallskip
After Cantor, Dedekind, Hausdorff, Bourbaki and up to the last decades we have been always moving in reverse direction.
Not anymore.

\smallskip

To help an uninitiated reader to see what is going on,
here are simple examples:  imagine that  ``two''  counts not fingers or
stones, but  orientations of an Euclidean space,
and $\bold{Z}$ counts homotopy classes of closed oriented loops around
zero in a plane. 

\smallskip
Notice that in the Set Theory, where we, after Cantor,
Frege and Russell, interpret integers as cardinalities of finite  discrete sets,
only natural numbers $\bold{N}$  appear directly,
whereas negative numbers require a psychologically 
difficult and historically late  leap of imagination.

\medskip

{\bf A.3. Renormalization of financial markets?} Finally, it might be amusing to think about the current financial crisis in
similar terms. 

\smallskip

This is not a simple fancy: recall that one of the early psychologically acceptable interpretations of negative numbers was formulated in terms 
of {\it debt.}

\smallskip

Now, to put it crudely, consider  the sum total of values of contracts defining what all
players at the global financial market {\it owe} to their counterparts 
(``minus infinity'') at a given time, and the respective
sum total of values of contracts defining what all
players at the financial market {\it are owed} (``plus infinity''). In the 
world of  material reality where one can lend only what one owns
and get back only what one has lent, these two ``infinities
for dummies''  would  exactly cancel.

\smallskip

Money adds to this material world an ideal dimension of {\it credit}, and (in good times)
a creative force. The difference between credits and debits 
pays for a private house and public education.
In a healthy economy, this difference however must 
be reasonably stable and remain on a considerably lesser scale
than the two infinities.

\smallskip

Technically, debts and credits
do not cancel anymore for many reasons. Banks are required to hold monetary reserves
which constitute a small percentage of their deposits; the rest can be invested, loaned etc.
 Debts must be paid 
at different times in   future, at various rates of interest.
This line of thought is iterated, which leads to the creation of contracts buying and selling risks, debts, etc. Such contracts are {\it derivatives},
financial instruments whose value is derived from the value of 
other financial instruments. The two infinities, and their difference,
start fluctuating and eventually lose their  contacts with reality.{\footnotemark1}
\footnotetext{Cf.  [Mas] for a related quantitative discussion.}

\smallskip

In a remarkable agreement with our metaphor, one of the great players and keen observers of financial markets,
Warren Buffett, once called derivatives ``financial weapons of mass destruction''.

\smallskip

Renormalization of finances probably
needs intellects of Richard Feynman's scale.

\medskip

{\bf A.4. Computational viewpoint on human civilization.} From the computational viewpoint,  human civilization is a supplier of software, hardware, and oracles
producing programs and inputs. It is also a consumer of outputs.

\smallskip

{\it Scientific laws} such as Newton's law of gravity, are 
short (Kolmogorov simple) oracular prescriptions for writing software
that will be
calculating (predicting), say, visible movements of planets.

\smallskip 
{\it Scientific observations} are systematic methods of obtaining
oracular prescriptions for inputs into resulting programs. 
\smallskip

Gods write equations of motion, devils choose initial and boundary
conditions, experimenters store them in databases. 

\bigskip

\centerline{\bf References}

\medskip

 [AlMar1]  P.~Aluffi, M.~Marcolli. {\it Feyman motives of banana graphs.} hep--th/0807.1690 
 \smallskip
 
 [AlMar2]  P.~Aluffi, M.~Marcolli. {\it Algebro--geometric Feynman rules.} hep--th/0811.2514
 
 \smallskip
 
 [AlMar3]  P.~Aluffi, M.~Marcolli. {\it Parametric Feynman integrals and determinant hypersurfaces.}
 math.AG/0901.2107
 
 \smallskip

[BeDr] A.~Beilinson, V.~Drinfeld. {\it Chiral Algebras.}  AMS Colloquium
Publications,  vol. 51, AMS, Providence, Rhode Island, 2004.
\smallskip
[Bl] S.~Bloch. {\it Motives associated to graphs.} Japan J.~Math., vol.~2 (2007), 165--196.

\smallskip
[BlEsKr] S.~Bloch, E.~Esnault, D.~Kreimer. {\it On motives associated to
graph polynomials.} Comm.~Math.~Phys., vol. ~267 (2006), 181--225.

\smallskip

[BoMa] D.~Borisov, Yu.~Manin.  {\it Generalized operads and their inner cohomomorhisms.}
 In: Geometry and Dynamics of Groups
and spaces (In memory of Aleksader Reznikov). Ed. by M. Kapranov et al.
Progress in Math., vol. 265. 
Birkh\"auser, Boston, pp. 247--308.
math.CT/0609748

[BdeM1] L.~Boutet de Monvel. {\it Remark on divergent multizeta series.}
Posted at http://people.math.jussieu.fr/$\sim$boutet, 2004.

\smallskip

[BdeM2] L.~Boutet de Monvel. {\it Alg\`ebre de Hopf des diagrammes de
Feynman, renormalisation et factorisation de Wiener--Hopf (d'apr\`es
A.~Connes et D.~Kreimer)}. S\'eminaire Bourbaki,
no. 900, 2002.

\smallskip

[BrYa] M.~Braverman, M.~Yampolsky. {\it Computability of Julia sets.} Moscow Math. Journ., 8:2 (2008), 185--231.

\smallskip

[Cas1] D.~Castella. {\it L'Alg\`ebre tropicale comme alg\`ebre de la caracteristique 1:
Alg\`ebre lin\'eaire sur les semi--corps idempotents.} math.AC/0807.3088

\smallskip

[Cas2] D.~Castella. {\it L'Alg\`ebre tropicale comme alg\`ebre de la caracteristique 1:
Polyn\^omes rationnels et fonctions polynomiales.} math.RA/0809.023

\smallskip

[ConKr1] A.~Connes, D.~Kreimer. {\it Hopf algebras, renormalization and noncommutative geometry.}
hep-th/9808042

\smallskip

[ConKr1] A.~Connes, D.~Kreimer. {\it Renormalization in quantum field  theory
and the Riemann--Hilbert problem. I. The Hopf algebra structure of
graphs and the main theorem.} Comm. Math. Phys. 210, no. 1 (2000),
249--273.

\smallskip

[Cos] K.~Costello. {\it Renormalization of quantum dield theories.}
Available at author's home page.
\smallskip

[De] P.~Deligne. {\it Multiz\'eta values.} Lecture notes, IAS Princeton, 2001.

\smallskip

[E-FMan] K.~Ebrahimi--Fard and D.~Manchon. {\it The combinatoris of
Bogolyubov's recursion in renormalization.} math-ph/0710.3675

\smallskip

%[Fr1]  A.~Frabetti. {\it Groupes de s\'eries et renormalisation des champs quantiques.}
%\smallskip

[Fra]  A.~Frabetti.  {\it Renormalization Hopf algebras and combinatorial groups.}  

arxiv:0805.4385

\smallskip

[Fri] B.~Friedrich. {\it Periods and algebraic de Rham cohomology.}
math.AG/0506113

\smallskip

[GaLe] P.~G\'acs, A.~Levin. {\it Causal Nets or What Is a Deterministic Computation?}
Int.~Journ.~Theor.~Phys., vol. 21, No. 12 (1982), 961--971.

\smallskip

[GeMa] S.~Gelfand, Yu.~Manin. {\it Methods of homological
algebra.} Second edition, Springer, 2003.

\smallskip

[Gl] M.~D.~Gladstone. {\it A reduction of the recursion scheme.}
Journ. Symb. Logic, 32:4 (1967), 505--508.

\smallskip
[I1] L.~Ionescu. {\it The Feynman legacy.} math.QA/0701069

\smallskip

[I2] L.~Ionescu. {\it From operads and PROPs to Feynman processes.}
math.QA/0701299

\smallskip

[Ka] M.~Kapranov. {\it Noncommutative geometry and path integrals.}
arXiv:math/0612411

\smallskip
[KoZa] M.~Kontsevich, D.~Zagier. {\it Periods.} In: Mathematics unlimited ---2001 and beyond, Springer, Berlin 2001, 771--808.

\smallskip

[Kr1] D.~Kreimer. {\it On the Hopf algebra structure of perturbative
quantum field theory.} Adv. Theor. Math. Phys., 2 (1998), 303--334.

\smallskip

[Kr2] D.~Kreimer. {\it Structures in Feynman graphs  -- Hopf algebras and
symmetries.}   hep-th/020211

\smallskip 

[Kr3] D.~Kreimer. {\it Dyson--Schwinger equations: from Hopf algebras to
number theory.} hep--th/0609004

\smallskip

[KremSz] K.~Kremnizer, M.~Szczesny. {\it Feynman graphs, rooted trees,
and Ringel--Hall algebras.}  math.QA/0806.1179

\smallskip

[Le] L.~Levin. {\it Various measures of complexity for finite objects
(axiomatic description).} Soviet Math. Dokl., vol 17, No. 2 (1976), 522 --526.

\smallskip

[Manch] D.~Manchon. {\it Hopf algebras, from basics to
application to renormalization.} math.QA/0408405

\smallskip

[Ma1]  Yu.~Manin.  {\it A Course in Mathematical Logic. } Springer Verlag, 1977. XIII+286 pp. (The second, expanded Edition in preparation). 

\smallskip

[Ma2]  Yu.~Manin. {\it Classical computing, quantum computing,
and Shor's factoring algorithm.}  S\'eminaire Bourbaki, no. 862 (June 1999),
Ast\'erisque, vol 266, 2000, 375--404.
 quant-ph/9903008.
 
 \smallskip
 
 [Ma3] Yu.~Manin. {\it Frobenius manifolds, quantum cohomology, and moduli
spaces.}  AMS Colloquium Publications, vol. 47, Providence, RI, 1999,
xiii+303 pp.

\smallskip

[Ma4] Yu.~Manin. {\it Renormalization and computation II:
Time cut--off and the Halting Problem. math.QA/0908.34.30}

\smallskip

[Mar] M.~Markl. {\it Operads and PROPs.} math.AT/0601129

\smallskip

[MarkShi] I.~L.~Markov, Y.~Shi. {\it Simulating quantum computation by contracting
tensor networks.}  quant--ph/0511069

\smallskip

[Mas] V.~Maslov. {\it Economic law of increase of Kolmogorov complexity.
Transition from financial crisis 2008 to the zero--order phase
transition (social explosion).} q-fin.GN/0812.4737
 
\smallskip

[Mi] R.~Miller. {\it Computable fields and Galois theory.} Notices of the AMS,
vol.~55, no.~7 (2008), 798--807.

\smallskip
 
 [Mo] Y.~N.~Moschovakis. {\it What is an Algorithm.} Available
 at author's home page.

\smallskip
[Po] M~Polyak. {\it Feynman diagrams for pedestrians and mathematicians.}

 math.GT/0406251
\smallskip

[PaSt1] L.~Pachter, B.~Sturmfels. {\it Tropical geometry of
statistical models.} Proc.~Nat. Ac.~Sci.~USA, vol. 101, no. 46 (2004),
16132--16137.

\smallskip

[PaSt2] L.~Pachter, B.~Sturmfels. {\it  Parametric inference for biological sequence
analysis.} q-bio.GN/0401033

\smallskip

[Ra] G.~Racinet. {\it Doubles m\'elanges des polylogarithmes multiples aux
racines del'unit\'e.}  Publ. Math. Inst. Hautes \'Etudes Sci.,  No. 95  (2002), 185--231. 

\smallskip

[Ro] R.~M.~Robinson. {\it Primitive recursive functions.}
Bull. AMS, 53 (1947), 925--942.

\smallskip

[Sc] D.~Scott. {\it The lattice of flow diagrams.} In: Symposium on Semantics of
Algorithmic Languages, Springer LN of Mathematics, 188 (1971), 311--372.

\smallskip

[TsVlN] M.~Tsfasman, S.~Vl\v{a}du\c{t}, D.~Nogin. {\it Algebraic geometric codes: basic notions.} Mathematical Surveys and Monographs, 139. American Mathematical Society, Providence, RI, 2007. xx+338 pp. 

\smallskip

[Va1] B.~Vallette. {\it A Koszul duality for PROPs.} math.AT/0411542

\smallskip

[Va2] B.~Vallette. {\it Manin products, Koszul duality, Loday algebras and Deligne conjecture.}
math.QA/0609002

\smallskip

[vdL] P.~van der Laan. {\it Operads and the Hopf algebras
of renormalization.} math-ph/0311013

\smallskip

[VlMa]  S.~Vl\v{a}du\c{t}, Yu.~ Manin. {\it Linear codes and modular curves. (Russian).} In:  Current problems in mathematics, Vol. 25,  209--257, Itogi Nauki i Tekhniki, Akad. Nauk SSSR,  
VINITI, Moscow, 1984.

\smallskip

[Ya] N.~S.~Yanofsky. {\it Towards a definition of an algorithm.}
math.LO/0602053

\smallskip

[Yo] M.~Yoshinaga. {\it Periods and elementary real numbers.}
math.AG/0805.0349

\smallskip

[Zo]  P.~Zograf. {\it Tensor networks and the enumeration of regular subgraphs.}
math.CO/0605256

\enddocument